\documentclass[amstex]{article}
\usepackage{mathptmx}
\usepackage{helvet}
\usepackage{courier}
\usepackage{graphicx}
\usepackage{amssymb}
\usepackage{subfig}
\usepackage{makeidx}
\usepackage[all]{xy}
\usepackage{multicol}
\pdfoutput=1
\makeindex             

\usepackage{footmisc}
\usepackage{epstopdf}
\usepackage{url}
\usepackage{multicol}
\newtheorem{theorem}{Theorem}[section]
\newtheorem{cor}[theorem]{Corollary}
\newtheorem{prop}[theorem]{Proposition}
\newtheorem{thm}[theorem]{Theorem}

\newtheorem{notation}[theorem]{Notation}
\newtheorem{example}[theorem]{Example}
\newtheorem{remark}[theorem]{Remark}
\newtheorem{lemma}[theorem]{Lemma}
\newtheorem{definition}[theorem]{Definition}

\newenvironment{defn}{\begin{definition}\rm}{\end{definition}}
\newenvironment{notn}{\begin{notation}\rm}{\end{notation}}
\newenvironment{ex}{\begin{example}\rm}{\end{example}}

\newcommand{\Z}{\mathbb{Z}}

\newcommand{\A}{\mathbb{A}}
\newcommand{\N}{\mathbb{N}}
\newcommand{\Q}{\mathbb{Q}}

\newcommand{\m}{\mathfrak{m}}

\newcommand{\n}{\mathfrak{n}}

\newcommand{\IC}{\mathbb{C}}

\newcommand{\relint}{\mathrm{relint}}

\newcommand{\conv}{\mathrm{conv}}
\newcommand{\ord}{\mathrm{ord}}
\newcommand{\spec}{\mathrm{Spec}}
\newcommand{\pic}{\mathrm{Pic}}

\usepackage[mathcal]{euscript}
\def\qedbox{\hbox{$\rlap{$\sqcap$}\sqcup$}}

\begin{document}
\title{Weak Normality and Seminormality}
\author{Marie A. Vitulli  \\
 Department of Mathematics \\
1222 University of Oregon \\
Eugene, OR 97403-1222} 
\date{June 17, 2009}                                           
\maketitle
 \abstract{In this survey article we outline the history of the twin theories of weak normality and seminormality for
 commutative rings and algebraic varieties with an emphasis on the recent developments in these theories
  over the past fifteen years.   We develop the theories for general commutative rings, but specialize to reduced Noetherian rings when necessary.  We hope to acquaint the reader with many of the consequences of the theories. }
  
  \medskip
  
   All rings in this paper are commutative with identity, all modules are unitary, and ring homomorphisms preserve the identity.

\section{Introduction}

The operation of \textbf{weak normalization} was formally introduced in 1967 by A. Andreotti and F. Norguet 
\cite{Andreotti:1967aa} in order  to solve a problem that arose while constructing a certain parameter space associated with a complex analytic variety.  Their construction was dependent on
the embedding of the  space in complex projective space.  In this setting the normalization of the parameter space is independent of the embedding, but no longer parametrizes what it was intended to since one point may split into several in the normalization.  To compensate one ``glues'' together points on the normalization that lie over a single point in the original space.  This leads to the weak normalization of the parameter space  whose underlying point set is in one-to-one correspondence with the point set of the parameter space.  
\index{weak normalization} 
 \index{Andreotti, Aldo!weak normalization} \index{Norguet, Fran\c{c}ois!weak normalization of complex  space}
A few years later weak normalization was introduced in the context of schemes and their morphisms by
 A. Andreotti and E. Bombieri.  For an integral extension of a local ring
  $A \subset B$,  they first introduced the notion of ``gluing" \index{gluing}the prime ideals of $B$ lying over the unique
   maximal ideal of $A$,  mirroring the complex analytic construction. This notion of gluing, which we will refer to as ``weak gluing,"  \index{weak gluing} appears in Lemma \ref{lemma on weak local gluing} and is formally defined in Definition \ref{defn weak gluing} . This algebraic operation results in a
    local ring, integral over the given local ring, and whose residue class field is purely inseparable over the given residue class field  \cite{Andreotti:1969aa}.
     \index{Bombieri, Enrico!weak normalization of scheme}  \index{Andreotti, Aldo!weak normalization of scheme}
  For a general integral extension $A \subset B$, an element $b \in B$ is in the weak normalization of $A$ in $B$ if and only if  for every point $x \in \spec (A)$, the image of $b$ is in the weak gluing of $A_x$ in $B_x$ over $x$.   
Andreotti and Bombieri then turned their attention to schemes and their structure sheaves.  They assumed they were working with preschemes, \index{prescheme} however, what was called a prescheme in those days (e.g. see Mumford's Red Book \cite{MR1748380}) is today called a scheme so we dispense with the prescheme label.   They defined the \textbf{weak subintegral closure } of the structure sheaf pointwise using the notion of gluing they already defined for local rings.  
 They next defined the \textbf{weak normalization} $\sigma \colon X^* \to X$ of a reduced algebraic scheme $X$ over an arbitrary field $K$ so that the scheme ${}^*X$ represents the weak subintegral closure of  the structure sheaf of $X$.   Andreotti and Bombieri established a universal mapping property  of this pair $(X^*, \sigma)$ which we discuss in Section \ref{subsec: geom props}.
 
 \index{weak subintegral closure}  \index{universal mapping property}

At about the same time C. Traverso \cite{Traverso:1970zv} \index{Traverso, Carlo!seminormalization}introduced the closely related notion of the
 \textbf{seminormalization} \index{seminormalization}of a commutative ring $A$ in an integral extension $B$.    Like in the
  Andreotti-Bombieri construction, given a local ring $A$ one glues  the prime ideals of $B$ lying over the
   unique maximal ideal of $A$ (i.e,  the maximal ideals of $B$) but this time in a way that results in a local
    ring with residue class field \textit{isomorphic} to that of $A$ (see Definition \ref{defn of gluing}). 
For an arbitrary integral extension of rings $A \subset B$, an element $b \in B$ is in the seminormalization
 ${}^+_BA$ of $A$ in $B$ if and only if for every point $x \in \spec (A)$, the image of $b$ is in the gluing of
  $A_x$ in $B_x$ over $x$.  Traverso showed \cite[Theorem 2.1]{Traverso:1970zv} that for a finite integral 
  extension of Noetherian rings $A \subset B$, the seminormalization ${}_B^+ A$  of $A$ in $B$ is obtained
   by a finite sequence of gluings.  Traverso defined a \textbf{seminormal ring} to be a ring that is equal to
    its seminormalization  in $\overline{A}$, where $\overline{A}$ denotes the integral closure of $A$ in  its total ring of quotients.  Traverso showed that for a
     reduced Noetherian ring $A$ with finite normalization and a finite number of indeterminates $T$, the
      ring $A$ is seminormal if and only if the canonical homomorphism $\pic (A) \to  \pic (A[T])$ is surjective. Here $\pic (A)$ denotes the \textbf{Picard group} of $A$, namely the group of isomorphism classes of rank one projective modules over $A$ with $\otimes_A$ as the group law. 
 \index{seminormal ring} \index{Traverso, Carlo!seminormal ring}
 \index{Picard group} \index{Traverso, Carlo!seminormality and Picard group} \index{$\pic (A)$: Picard group of $A$} 
 
A few years after Traverso's paper appeared E. Hamann \index{Hamann, Eloise} showed that a seminormal ring $A$ contains
 each element $a$ of its total quotient ring such that $a^n, a^{n+1} \in A$ for some positive integer $n$
  (see \cite[Prop. 2.10]{Hamann:1975hq}). \index{Hamann, Eloise!Hamann's criterion for seminormality}Hamann showed that this property is a characterization of
   seminormality for a pseudogeometric ring ( \cite[Prop. 2.11]{Hamann:1975hq}); today pseudogeometric rings are more commonly known as Nagata rings.   More generally,
    Hamann's criterion characterizes seminormal reduced Noetherian rings.  It also can be used to
     characterize arbitrary seminormal extensions, as we shall see. The integers $n, n+1$ that appear in
      Hamann's characterization can be replaced by any pair $e, f$ of relatively prime positive integers
       (\cite[Prop. 1.4]{Leahy:1981ta}.

A slight but significant refinement of Hamann's criterion was made by R. Swan  \cite{Swan:1980kn}  in 1980.  Swan defined a \textbf{seminormal ring} \index{seminormal ring} as a reduced ring $A$ such that whenever $b, c \in A$ satisfy 
$b^3 = c^2$ there exists $a \in A$ such that $b= a^2, c = a^3$.  For reduced Noetherian rings (or, more generally, any reduced ring whose total quotient ring is a product of fields), Swan's definition of a seminormal
 ring is equivalent to Traverso's.  It was pointed out by D.L. Costa \cite{Costa:1982a} \index{Costa, Douglas L.} that a ring that satisfies
  Swan's criterion is necessarily reduced so the $b^3 = c^2$ criterion alone characterizes seminormality. Swan 
  showed in \cite[Theorem 1]{Swan:1980kn} that $\pic  (A) \cong \pic  (A[T])$ for some finite set of indeterminates
   $T$ is equivalent to the seminormality of $A_{\mathrm{red}}$.  In \cite{MR558489} Gilmer and Heitmann     presented an example of a non-Noetherian reduced ring that is equal to its own total ring of quotients but such that $\pic (A)$
     is not canonically isomorphic to $\pic (A[T])$. Thus according to Traverso's original definition of a seminormal ring, the seminormality of $A$ isn't equivalent $\pic \, A \cong \pic \,  A[T]$.   In  \cite{Swan:1980kn} Swan constructed the seminormalization of a general commutative ring in a way that is reminiscent of the
  construction of the algebraic closure of a field.  He went on to show that any reduced commutative ring has an essentially unique seminormalization  (a subintegral extension $A \subset {}^+A$ such that ${}^+A$  is a seminormal ring). 
  \index{Swan, Richard!seminormal ring}  
\index{Gilmer, Robert} \index{Heitmann, Raymond}

Early in the 1980s, both Greco-Traverso \cite{MR571055} \index{Greco, Silvio} \index{Traverso, Carlo} and Leahy-Vitulli \cite{Leahy:1981ta} \index{Leahy, John V.} \index{Vitulli, Marie A.} published pivotal papers that linked the earlier work of Andreotti-Norguet-Bombieri to the work of Traverso and looked at
 the singularities of schemes and varieties.  Both pairs of authors showed that a reduced, complex analytic
  space is weakly normal at a point $x \in X$ in the sense of Andreotti-Norguet  if and only if the local ring of 
  germs of holomorphic functions $\mathcal{O}_{X,x}$ is seminormal  in the sense of Traverso 
  (cf. \cite[Prop. 2.24]{Leahy:1981ta} and \cite[Cor. 5.3]{MR571055}).   Leahy-Vitulli   defined  a  weakly normal singularity  called a multicross \index{multicross}
    (see \cite{Leahy:1981ta} and \cite{Leahy:1981kr}). \index{Leahy, John V.!multicross} \index{Vitulli, Marie A.!multicross} Briefly a point $x$ on a variety $X$ is a \textbf{multicross} if $x \in X$ is
      analytically isomorphic to $z \in Z$ where $Z$ is the union of linearly disjoint linear subspaces.  Most singularities of a weakly normal variety are of this
     type in the sense that the complement of the set of multicrosses forms a closed subvariety of codimension at least 2 
     \cite[Theorem 3.8 ]{Leahy:1981kr}.   Leahy and Vitulli worked with  algebraic
   varieties over an algebraically closed field of characteristic 0.
\index{Andreotti, Aldo} \index{Norguet, Fran\c{c}ois}

A few years later, Yanagihara  gave an intrinsic definition of a \textbf{weakly normal ring} \index{weakly normal ring} analogous to the Traverso-Swan definition of a seminormal ring (see Definition \ref{Yanagihara intrinsic def wn} or  \cite{Yanagihara:1985os} ).  Yanagihara said that a reduced ring $A$ is
 weakly normal provided that the ring is seminormal in the sense of Swan and another condition involving rational primes holds.
\index{Yanagihara, Hiroshi!weakly normal ring}

Before we can talk about the next set of results we need to recall a pair of definitions.  An integral extension of rings $A \subset B$ is called a \textbf{(weakly) subintegral extension} if for each prime ideal $P$ of $A$ there is a unique prime ideal $Q$ of $B$ lying over $P$ and the induced map
  of residue fields $A_P / PA_P \to B_Q/ QB_Q$ is an isomorphism (a purely inseparable extension) (see Definitions \ref{defn subintegral ext} and \ref{defn weakly subintegral ext}).  We point out the (weak) subintegral closure plays the role in (weak normality) seminormality that integral closure plays in normality.
   \index{subintegral extension} \index{weakly subintegral extension}
   
  In the mid to late 1990s a series of papers by L. Reid, L. Roberts, and B. Singh appeared in which an elementwise criterion for weak subintegrality was introduced and developed.  At first they worked with $\Q$-algebras and introduced what they called a system of subintegrality (SOSI).  They proved that an element $b$ in an integral extension $A \subset B$  admits a SOSI if and only if $b$ is subintegral over $A$.  A system of subintegrality is global in nature and makes no reference to local gluings.  We view the SOSI  as the first true ``elementwise"  criterion for subintegrality over a $\Q$-algebra.   Over the years Reid, Roberts, and Singh showed each of many conditions is equivalent to the existence of a system of subintegrality.  Among the equivalent conditions is  the existence of a highly-structured sequence of monic polynomials  of various degrees such that the element $b$ is a common root of each of these polynomials.   
\index{system of subintegrality (SOSI)} \index{Roberts, Leslie G.!SOSI} \index{Singh, Balwant!SOSI} \index{Reid, Les!SOSI} \index{Roberts, Leslie G.!criterion for weak subintegrality} \index{Singh, Balwant!criterion for weak subintegrality} \index{Reid, Les!criterion for weak subintegrality}

Reid, Roberts, and Singh were eventually able to drop the assumption that the rings are $\Q$-algebras, but only if they talked about weak subintegrality rather than subintegrality.   However they still talked about SOSIs, which perhaps should have been renamed systems of \textit{weak}  subintegrality.

In 1999 Roberts observed that the lower degree polynomials are rational multiples of the polynomial of highest degree.  Gaffney and Vitulli \cite[Prop. 2.2 and Prop. 2.3]{GaffVit2008:aa} \index{Gaffney, Terence} \index{Vitulli, Marie A.} have recently expanded this observation to make some geometric sense out of the sequence of equations that appears in the work of Reid, Roberts, and Singh.  
Proposition \ref{nec for wn} in this paper is a new result that gives another elementwise characterization of  weak subintegrality. \index{Vitulli, Marie A.! criterion for weak subintegrality} This also provides a simple, purely algebraic explanation, for why an element satisfying the sequence of equations introduced by Roberts, Reid, and Singh is necessarily weakly subintegral over the base ring.

The elementwise criteria of Reid, Roberts, and Singh enabled the current author and Leahy to talk about when an element $b$ of a ring $B$  is weakly subintegral over an ideal  $I$ of a subring $A$ and define the weak subintegral closure of an ideal either in the containing ring $A$ or an integral extension ring $B$ \cite{Vitulli:1999fi}.   Vitulli and Leahy showed that for an extension $A \subset B$ of rings, an ideal $I$ of $A$, and $b \in B$, the element $b$ is weakly subintegral over $I^m$ if and only if $bt^m$ is weakly subintegral over the Rees ring $A[It]$ \cite[Lemma 3.2]{Vitulli:1999fi}.  Thus the weak subintegral closure ${}^*_B I$ of $I$ \index{weak subintegral closure of an ideal}  in $B$ is an  
ideal of the weak subintegral closure ${}^* _B A$ of $A$ in $B$ \cite[Prop. 2.11]{Vitulli:1999fi}.   In particular, the weak 
subintegral closure of $I$ in $A$ is again an ideal of $A$, which we denote simply by ${}^*I$ .  Vitulli and Leahy showed that given an ideal $I$ in a reduced ring $A$ with finitely many minimal primes and total quotient ring $Q$, it holds that ${}^*(A[It]) = \oplus_{n \ge 0} {}^* _Q (I^n)t^n$ \cite[Cor. 3.5]{Vitulli:1999fi}. Gaffney and  Vitulli \cite{GaffVit2008:aa} further developed the theory of weakly subintegrally closed ideals in both the algebraic and complex analytic settings. They defined a subideal $I_>$ of the weak subintegral closure ${}^*I$ of an $I$ in a Noetherian ring that can be described solely in terms of the Rees valuations of the ideal and related this subideal to the minimal reductions of $I$.  Gaffney and  Vitulli also proved a valuative criterion for when an element is in the weak subintegral closure of an ideal.  \index{Gaffney, Terence} \index{Vitulli, Marie A.} H. Brenner \cite{Brenner2006aa} \index{Brenner, Holger} has proposed another valuative criterion in terms of maps into the appropriate test rings. 

We now describe the contents of this paper.   Section 2 deals with seminormality and seminormalization.  In Section \ref {subsect:si and sn ext and rel sn} we recall Traverso's notion of gluing the maximal ideals of an integral extension of a local ring and prove the fundamental properties that the ring so obtained 
enjoys.   We give Traverso's definition of the seminormalization of a ring $A$ in an integral extension $B$
 of $A$ and provide proofs of the fundamental properties of seminormalization.  In Section
  \ref{subsect:swan's refinement} we take a look at Swan's refined notion of a seminormal ring and
 his construction of the  seminormalization of a ring.    We discuss seminormality and Chinese Remainder Theorems in Section \ref{subsec:CRT}.   In Section \ref{subsect:sn and picard grp} we explore the connections between
seminormality and when rank one projectives over the polynomial ring $A[X]$ are extended from rank
one projectives over $A$.   In Section \ref{subsect:1 diml sn rings} we recall some characterizations of
 seminormal one-dimensional local rings.    In Section \ref{sect:weak normality and weak normalization} we develop the theories of  weakly subintegral extensions and  weak normality.  We discuss weakly subintegral and weakly normal extensions and the operation of weak normalization in Section \ref{subsect:wsi and wn extensions}.  Section \ref{subsect:SOSI} looks at the systems of subintegrality introduced by Roberts and Singh and developed by Reid, Roberts, and Singh.  We offer a new algebraic criterion for the weak subintegrality of an element over a subring in Section \ref{subsect:fresh approach}.  We give the reader a very brief glimpse of some geometric aspects of weakly normal varieties in Section \ref{subsec: geom props}.   In Section \ref{subsect:wn and CRT} we recall a couple of Chinese Remainder Theorem results for weakly normal varieties.   Our final section, Section \ref{subsec: wsi of ideal} is devoted to the notions of weak subintegrality over an ideal and the weak subintegral closure of an ideal.  We recall the Reid-Vitulli geometric and algebraic characterizations of the weak subintegral closure of a monomial ideal.  We introduce the ideal $I_>$ that was defined by Gaffney and Vitulli  and cite some of their results connected with this ideal.  

\section{Seminormality and Seminormalization}

In Section \ref{subsect:si and sn ext and rel sn} we will outline the development  of gluings of prime ideals, subintegral extensions, and the seminormalization of a ring in an integral extension ring, a relative notion dependent on the extension ring. In Section \ref{subsect:swan's refinement} we will deal with the absolute notions of a seminormal ring and seminormalization; these  notions that do not depend on particular integral extensions.  In Section  \ref{subsec:CRT} we discuss seminormality in relation to various versions of the Chinese Remainder Theorem.  In the Section  \ref{subsect:sn and picard grp}   we discuss Picard group results connected with seminormality, which some think of as arithmetic consequences of the theory.  Most of the original proofs of these results were $K$-theoretic in nature.  In the Section \ref{subsect:1 diml sn rings} we discuss some results on one-dimensional seminormal rings.

We work with arbitrary commutative rings in this section.  We do not assume the rings we are discussing are reduced, Noetherian, or have any other special properties.

\subsection{Subintegral and Seminormal Extensions and Seminormalization Relative to an Extension} \label{subsect:si and sn ext and rel sn}
We start by giving the local construction that leads to the seminormalization of a ring in an integral extension.
For a ring $B$ we let $R(B)$ denote the Jacobson radical of $B$.  We notationally distinguish between a prime ideal $P = P_x$ of $B$ and the corresponding point $x \in \spec (B)$, when this is convenient.  For an element $b \in B$ and $x \in \spec (B)$ \index{$\spec (A)$: prime spectrum of $A$} we let $b(x)$ \index{$b(x)$: image of $b$ in $\kappa(x)$} denote the image of $b$ in the residue class field $\kappa(x) = B_x/ P_x B_x$. \index{$\kappa(x)$: residue field of $A_x$}  If $A \subset B$ are rings and $x \in \spec(A)$ we let $B_x$ denote the ring obtained by localizing the $A$-module $B$ at the prime ideal $x$. \index{$R(B)$: Jacobson radical of $B$}

\begin{lemma} \label{local gluing} Let $(A , \m)$ be a local ring and $A \subset B$ be an integral extension. 
Set $A' = A + R(B) \subset B$ and $\m ' = R(B)$.  For each maximal ideal $\n _i$ of $B$ let $\omega _i \colon A/ \m \to B/ \n_i$ be the canonical homomorphism.  The following assertions hold.
\begin{enumerate}
\item  $(A' , \m ' )$ is a local ring and the canonical homomorphism $A/ \m \to A'/\m '$ is an isomorphism;
\item  $A'$ is the largest intermediate  local ring $(C, \n)$ such that  $A/ \m \cong C / \n$; and
\item An element $b \in B$ is in $A'$ if only if 
\begin{enumerate}
\item$ b(x_i) \in \omega_i(\kappa(x)) \mbox{ for all closed points } x_i $ of $\spec(B)$, and
\item $\omega _i^{-1} (b(x_i)) = \omega _j^{-1} (b(x_j)) \mbox{ for all } i, j $.
\end{enumerate}
\end{enumerate}
\end{lemma}

\noindent \textbf{Proof.}

  1.  We first show that $A' = A + R(B)$ is local with unique maximal ideal $\m ' =  R(B)$.  It is clear that the canonical map  $A/ \m \to A' / \m '$ is an isomorphism and hence $\m '$  is a maximal ideal of $A'$.  Suppose that $\n$ is any maximal ideal of $A'$.  By Lying Over there exists a necessarily maximal ideal $\n_i$ of $B$ lying over $\n$.  Since $\n = \n _i \cap A' \supset R(B) \cap A' = R(B) = \m '$ we may conclude that $\m '$ is the unique maximal ideal of $A'$.

2. Suppose that $A \subset C \subset B$ are rings and that $(C, \n)$ is local with residue class field isomorphic to that of $A$.  As in the proof of 1. we must have $\n \subset R(B)$ and since $A/ \m \to C/ \n$ is an isomorphism we may conclude that $C \subset A + R(B)$.  

3.  Suppose that $b \in B$ satisfies a) and b) above.  Choose $a \in A$ such that $\omega_i(a(x)) = b(x_i)$ for some $i$, and hence for every $i$.  Then $b - a \in \n_i$ for every maximal ideal $\n _i$ of $B$, and hence,
$b \in A + R(B) = A'$.  It is clear that if $b \in A'$ then $b$ satisfies conditions  a. and b. above. \qedbox

\begin{definition} \label{defn of gluing} Let $(A , \m)$ be a local ring and $A \subset B$ be an integral extension.  We say that $A + R(B)$ is the \textbf{ring obtained from \emph{A} by gluing} the maximal ideals in $B$ over $\m$ or that $A + R(B)$ is the ring obtained from $A$ by gluing in $B$ over $\m$. Letting $x \in \spec(A)$ denote the point  corresponding to the maximal ideal $\m$, this ring is sometimes denoted by ${}^+_x A$. 
\end{definition}
\index{ ${}^+_x A$: gluing of $A$ over $x$}   \index{gluing}
We can also glue the prime ideals of an integral extension ring lying over an arbitrary prime ideal.  We state the result and leave the proof up to the reader.

\begin{lemma} \label{prime gluing}  Let $A \subset B$ be an integral extension of rings and $x \in \spec(A)$.  For each point $x_i \in \spec(B)$ lying over $x$, let $\omega _i \colon \kappa(x) \to \kappa(x_i)$ denote the canonical homomorphism of residue class fields.   Set $A' = \{ b \in B \mid b_x \in A_x + R(B_x) \}$.
The following assertions hold.
\begin{enumerate}
\item There is exactly one point $x' \in \spec(A')$ lying over $x$ and the canonical homomorphism $\kappa(x) \to \kappa(x')$ is an isomorphism. 
\item $A'$ is the largest intermediate ring with one point lying over $x$ and with isomorphic residue field at the corresponding prime; and
\item An element $b \in B$ is in $A'$ if and only if
\begin{enumerate}
\item$ b(x_i) \in \omega_i(\kappa(x)) \mbox{ for all  points } x_i $ of $\spec(B)$ lying over $x$, and
\item $\omega _i^{-1} (b(x_i)) = \omega _j^{-1} (b(x_j)) \mbox{ for all } i, j $.
\end{enumerate}
\end{enumerate}
\end{lemma}

We will now recall Traverso's global definition and its first properties.
\index{seminormalization}  
\begin{definition}  Let $A \subset B$ be an integral extension of rings.  We define the \textbf{seminormalization} \index{seminormalization} ${}^+_B A$ of $A$ in $B$ to be
\begin{equation}
{}^+_B A = \{ b \in B \mid b_x \in A_x + R(B_x) \mbox{ for all } x \in \spec(A) \}.
\end{equation}
\index{${}_B^+A$: seminormalization of $A$ in $B$}
\end{definition}
Before stating the first properties of ${}^+_BA$,  it will be useful at this point to recall  the notion of a subintegral extension.  Such extensions were earlier studied by Greco and Traverso in \cite{MR571055}, \index{Greco, Silvio} \index{Traverso, Carlo} where they were called quasi-isomorphisms,  but later on R. Swan \cite{Swan:1980kn} \index{Swan, Richard} called them subintegral extensions. 
\index{quasi-isomorphism}

\begin{definition} \label{defn subintegral ext}  A  \textbf{subintegral extension of rings } is an integral extension $A \subset B$ such that the associated map $\spec (B) \to \spec (A)$ is a bijection and induces isomorphisms on the residue class fields. An element $b \in B$ is said to be \textbf{subintegral} over $A$ provided that $A \subset A[b]$ is a subintegral extension.
\end{definition}
\index{subintegral extension}   \index{subintegral element} 

Notice that there are no proper subintegral extensions of fields.  We now look at some first examples of subintegral extensions.

\begin{ex}
Let $A \subset B$ be an extension of rings, $b \in B$ and $b^2, b^3 \in A$, then $A \subset A[b]$ is easily seen to be a subintegral extension.   If $A$ is any ring and $b, c \in A$ satisfy $b^3 = c^2$ then the extension $A \subset A[x] := A[X]/(X^2 - b, X^3 - c)$, where $X$ is an indeterminate, is thus a subintegral extension.
 \end{ex}
 
 \begin{defn} An   an \textbf{elementary subintegral extension} is a simple integral extension $A \subset A[b]$, such that $b^2, b^3 \in A$.
 \end{defn}
 \index{subintegral extension!elementary subintegral extension}
For an integral extension $A \subset B$ the seminormalization ${}^+ _ B A$ is the filtered union of all subrings of $B$ that can be obtained from $A$ by a finite sequence of elementary subintegral extensions.  Here by filtered union we mean that given any rings $C_1, C_2$ with the property there is a third ring $C$ with the property and satisfying $A \subset C_i \subset C \subset B  \;  (i=1,2)$. In particular the union of all such rings with the property is again a ring.
 
 \medskip
 The notion of a seminormal extension, which we now define, is \emph{complementary} to that of a subintegral extension.

\begin{definition} If $A \subset B$ is an integral extension of rings, we say $A$ is \textbf{seminormal }in $B$ if there is no subextension $A \subset C \subset B$ with $C \ne A$ and $A \subset C$ subintegral.

\end{definition}
\index{seminormal extension}

We point out that $A$ is seminormal in $B$ if and only if $A = {}^+_B A$.  We now recall Traverso's results regarding the seminormalization of $A$ in an integral extension $B$.

\begin{theorem}\cite[(1.1)]{Traverso:1970zv} \label{Traverso char} Let $A \subset B$ be an integral extension of rings. The following assertions hold.
\begin{enumerate}
 \item The extension $A \subset {}^+_B A$ is  subintegral;
 \item If $A \subset C \subset B$ and $A \subset C$ is subintegral, then $C \subset {}^+_B A$; 
 \item  The extension ${}^+_B A \subset B$ is seminormal; and
 \item  ${}^+_B A$ has no proper subrings containing $A$ and seminormal in $B$.
 
\end{enumerate}
\end{theorem}

\noindent \textbf{Proof.} 
1.  This follows immediately from Lemma \ref{prime gluing}.

2.  Now suppose that $A \subset C \subset B$ and that $A \subset C$ is subintegral.  Let $P \in \spec (A)$ and let $Q \in \spec (C)$ be the unique prime lying over $P$.  Since $A_P/ PA_P \to C_P/QC_P$ is an isomorphism and $QC_P \subset R(B_P)$ we must have $C_P \subset A_P + R(B_P)$. Since $P$ was arbitrary, $C \subset A'$.

3. This follows immediately from parts  1. and 2. 

4.  Suppose that $A \subset C \subset {}^+_B A$ and that $C \subset B$ is seminormal. Observe that $C \subset {}^+_B A$ is necessarily subintegral and hence $C = {}^+_B A$. \qedbox

\medskip

Paraphrasing this result, ${}^+_B A$ is the unique largest subintegral extension of $A$ in $B$ and is  minimal among the intermediate rings $C$ such that $C \subset B$ is seminormal. 

We mention some fundamental properties of  seminormal extensions.

\begin{prop} \label{prop conductors} Let $A \subset B \subset C$ be  integral extensions of rings.
\begin{enumerate}
\item  If $A \subset B$ is a seminormal extension, then the conductor $A:B$ of $A$ in $B$ is a radical ideal of $B$ (contained in $A$). 
\item   $A \subset B$ is seminormal if and only if $A:A[b]$ is a radical ideal of $A[b]$ for every $b \in B$.
 \item  If $A \subset B$ and $B \subset C$ are seminormal extensions, then so is $A \subset C$.
\end{enumerate}
\end{prop}
\noindent \textbf{Proof.}
For the proof of the first statement see  \cite[Lemma 1.3]{Traverso:1970zv} and for the second and third see
\cite[Prop. 1.4 and Cor. 1.5]{Leahy:1981ta}. \qedbox

\medskip

Before recalling a criterion due to Hamann and its generalizations we recall some definitions. 

\begin{definition}  Let $A \subset B$ be an integral extension of rings and $m, n$ be positive integers.  We say that $A$ is \textbf{\emph{n}-closed} in $B$ if $A$ contains each element $b \in B$ such that $b^n \in A$.  We say that $A$ is \textbf{(\emph{m,n})-closed} in $B$ if $A$ contains each element $b \in B$ such that $b^m, b^n \in A$.
\end{definition}
\index{$n$-closed}
\index{$(m,n)$-closed}
Hamann \cite[Prop. 2.10]{Hamann:1975hq} showed that a seminormal ring $A$ with integral closure $\overline{A}$ in its total quotient ring is $(n,n+1)$-closed in $\overline{A}$.   She also showed that this property is a characterization of seminormality for  pseudogeometric rings (also known as  Nagata rings) ( \cite[Prop. 2.11]{Hamann:1975hq}). We will that a ring $A$ (respectively, an integral extension $A \subset B$) satisfies  Hamann's criterion provided that $A$ is  $(n,n+1)$-closed in $\overline{A}$ (respectively, $B$), for some positive integer $n$.  Seminormality arose in Hamann's study of $R$-invariant and steadfast rings.  An $R$-algebra $A$ is \textbf{\textit{R}}-\textbf{invariant }provided that whenever $A[x_1, \ldots, x_n] \cong_R B[y_1, \ldots, y_n]$ for  indeterminates $x_i, y_j$ we must have $A \cong_R B$.  Finally, $R$ is called \textbf{steadfast} if the polynomial ring in one variable $R[x]$ is $R$-invariant. Hamann used Traverso's work on seminormality and the Picard group to show that a seminormal pseudogeometric local reduced ring is steadfast ( \cite[Theorem 2.4]{Hamann:1975hq}) and a seminormal domain is steadfast (\cite[Theorem  2.5]{Hamann:1975hq}).  
\index{$R$-invariant ring}  \index{steadfast ring}
Hamann's criterion characterizes reduced Noetherian rings that are seminormal in $\overline{A}$, where  $\overline{A}$ is as above. It also can be used to characterize arbitrary seminormal extensions as we shall see. The integers $n, n+1$ that appear in Hamann's characterization can be replaced by any pair $e, f$ of relatively prime positive integers (\cite[Prop. 1.4]{Leahy:1981ta}.

The following was proven by Leahy and Vitulli.

\begin{prop}  \cite[Prop. 1.4]{Leahy:1981ta} Let $A \subset B$ be an integral extension of rings.  The following are equivalent.
\begin{enumerate}
\item $A$ is seminormal in $B$;
\item $A$ is $(n,n+1)$-closed in $B$ for some positive integer $n$; and
\item $A$ is $(m,n)$-closed in $B$ for some relatively prime positive integers $m,n$. 
\end{enumerate}
\end{prop}

The reader should note if  $A$ is seminormal in $B$, then $A$ is $(m,n)$-closed in $B$ for \emph{every} pair of relatively prime positive integers $m$ and $n$. 
Recall from the introduction Hamann's results that assert  that for  any positive integer $n$ a seminormal ring is $(n,n+1)$-closed and that a pseudogeometric ring that is  $(n,n+1)$-closed is necessarily seminormal. 

The most commonly cited version of 2. above is with $n=2$, 	that is, the integral extension $A \subset B$ is seminormal if and only if it is $(2,3)$-closed (see \cite[Props. 2.10 and 2.11]{Hamann:1975hq}.  We now give our first  example of a seminormal extension.
 
 \begin{ex} \label{sn ext ex} Let $K$ be a field, $x$ be an indeterminate,  and consider $A =K[x^2] \subset B=K[x]$. We wish to see that $A$ is $(2,3)$-closed in $B$.  Let  $f = x^m(a_0+a_1x+a_2x^2 + \cdots +a_nx^n)$, where  $f^2, f^3 \in A$ and show that $f \in A$.  It suffices to assume that $a_0 = 1$ Consider $B = \oplus Kx^d$ as an $\N$-graded ring.  Then, $A$ is a graded subring.  Now, 
 \begin{eqnarray*}
f^2 &= & x^{2m} \left(1 + \sum_{k=1}^{2n} \sum_{i+j=k} a_ia_j x^k \right) \\
f^3 &=& x^{3m} \left(1 + \sum_{\ell = 1}^{3n} \sum_{i+j+k=\ell} a_ia_ja_k x^{\ell} \right).
 \end{eqnarray*}
 Looking at the lowest degree term of $f^3$ we see that $x^{3m} \in A$  and hence $m=2m'$ is even.  Observe that   $x^{4m'} \left(1 + \sum_{k=1}^{2n} \sum_{i+j=k} a_ia_j x^k \right) \in A \Rightarrow$  $ \sum_{k=1}^{2n} \sum_{i+j=k} a_ia_j x^k \in A \Rightarrow \sum_{i+j=k} a_ia_j =0$ for all odd numbers $k$ between 1 and $n$.
 
 First suppose char($K) \ne 2$.  Then, $2a_1 = 0 \Rightarrow a_1= 0$.  By induction on $k$, we may conclude $a_k = 0$ or all odd numbers $k$ between 1 and $n$.   Hence $f \in A$.
 
 Now suppose that char($K) = 2$.  Since $f^3 \in A$ we may conclude that \\
 $ \sum_{\ell = 1}^{3n} \sum_{i+j+k=\ell} a_ia_ja_k x^{\ell}  \in A$, which implies, $\sum_{i+j+k=\ell} a_ia_ja_k = 0$ for all odd numbers $k$ between 1 and $n$.  Starting with $\ell = 1$ we see that $3a_1 =a_1= 0$.   By induction, $a_k = 0$ or all odd numbers $k$ between 1 and $n$. Thus $f \in A$.  Therefore, $A$ is (2,3)-closed in $B$ and we may conclude that $A \subset B$ is a seminormal extension by Hamann's criterion.
 \end{ex}
With the previous characterization of a seminormal extension, it is easy to show that seminormal extensions are preserved under localization.

\begin{prop}\cite[Cor. 1.6]{Leahy:1981ta}  If $A \subset B$ is a seminormal extension and $S$ is any multiplicative subset, then $S^{-1}A \subset S^{-1}B$ is again seminormal.  Moreover, the operations of seminormalization and localization commute.
\end{prop}

Swan gave another useful characterization of a subintegral extension.
\index{Swan, Richard}

\begin{lemma} \cite[Lemma 2.1]{Swan:1980kn} \label{field char of subintegral}
An extension $A \subset B$ is subintegral if and only if $B$ is integral over $A$ and for all homomorphisms $\varphi \colon A \to F$ into a field $F$, there exists a unique extension $\psi \colon B \to F$.
\end{lemma}
\noindent \textbf{Proof.}  
 This follows from noting that $\varphi$ is determined by specifying a prime $Q$ (which will serve as  $\ker \psi)$ lying over $P = \ker \varphi$ and an extension of $\kappa(P) \to F$ to a map $\kappa(Q) \to F$. \qedbox

\subsection{Seminormal Rings, Seminormalization, and Swan's Refinements} \label{subsect:swan's refinement}

Traverso's original definition of a seminormal extension of rings and his  construction of the seminormalization of a ring $A$ in an integral extension $B$ of $A$ are still accepted by specialists in the field today. Since this conceptualization of seminormalization depends on the extension ring we'll refer to it as the relative notion. However, the current  \emph{absolute} notion of a seminormal ring and the construction of the seminormalization of a ring have changed due to Swan's insights. As the reader shall see, a seminormal ring in the sense of Swan is defined without any mention of an extension ring.  To say it differently, the relative theory is still due to Traverso but the absolute theory commonly used today is due to Swan.  We now recall from the introduction Traverso's original definition of a seminormal ring and Swan's modification of that definition.  One compelling reason for modifying the definition is that the modification enabled Swan to prove that  $A_{\mathrm{red}}$ is seminormal if and only if every rank one projective module over the polynomial ring $A[X]$ is extended from $A$ \cite[Theorem 1]{Swan:1980kn}. The reader is referred to Swan's original paper \cite{Swan:1980kn} for more details.  Lemma \ref{swan sn = traverso sn} below illustrates that under mild assumptions, Swan's absolute notion of a seminormal ring is equivalent to the earlier notion of Traverso.

For a ring $A$ we let $\overline{A}$ denote the integral closure of $A$ in its total ring of quotients.   
\index{ $\overline{A}$: integral closure of $A$ in its quotient ring}

\begin{definition} A \textbf{seminormal ring in the sense of Traverso} is a ring $A$ that is equal to its seminormalization in $\overline{A}$.
\end{definition}
\index{seminormal ring!seminormal ring in sense of Traverso}
Notice that if $A$ is a ring and $b$ is an element of the total ring of quotients of $A$ such that $b^2, b^3 \in A$, then $b$ is necessarily integral over $A$.  Thus a ring $A$ is seminormal in the sense of Traverso if and only if it is (2,3)-closed in its total quotient ring.

\begin{definition} A \textbf{seminormal ring in the sense of Swan} is a ring $A$ such that whenever $b, c \in A$ satisfy $b^3 = c^2$ there exists $a \in A$ such that $b= a^2, c = a^3$.  
\end{definition}
\index{seminormal ring!seminormal ring in the sense of Swan}
D. L. Costa \cite{Costa:1982a} \index{Costa, Douglas L.}observed that if $A$ is seminormal in the sense of Swan, then $A$ is necessarily reduced.  Just suppose $A$ is not reduced. Then, there is some nonzero element $b \in A$ such that $b^2 = 0$.  Then, $ b^3 = b^2 $ and hence there exists $a \in A$ such that $b = a^2, b = a^3$.  Hence $b = a^3= aa^2 = ab =aa^3 = b^2= 0$, a contradiction.  Swan's original definition of a seminormal ring stipulated that the ring is reduced, but in light of the previous remark this is redundant and we have omitted it from our definition.

Suppose that $A$ is seminormal in the sense of Swan and let $K$ denote its total ring of quotients. Assume that $b \in K$ and $b^2, b^3 \in A$.  Then, there exists an $a \in A$ such that $a^2 = b^2, a^3 = b^3$.  Then,  $(a - b)^3 = a^3 -3a^2b+3ab^2 - b^3 = a^3 -3b^3 + 3a^3 - b^3 = 0$ and hence $b = a \in A$, since $A$ is reduced. Thus $A$ must be seminormal in the sense of Traverso.  This same argument shows that for a seminormal ring $A$ and elements $b,c \in A$ satisfying $b^3 = c^2$, there exists a unique element $a \in A$ such that $b = a^2, c = a^3$.  The converse holds if $K$ is a product of fields, as we will now show.  

\begin{lemma} \label{swan sn = traverso sn}  Suppose that $A$ is a ring whose total ring of quotients $K$ is a products of fields.  Then, $A$ is seminormal in the sense of Traverso if and only if $A$ is seminormal in the sense of Swan.
\end{lemma}

\noindent \textbf{Proof.} Let $K = \prod K_i$, where $K_i$ are fields, and let $\rho_i \colon A \to K_i$ be the inclusion followed by the projection.   Notice that $A$ is reduced.  Assume that $A$ is seminormal in the sense of Traverso.  Suppose that $b, c \in A$ and $b^3 = c^2$.   If $b = 0$ then $b = 0^2, c= 0^3$.  Now assume $b \ne 0$.  Consider the element $\alpha = (\alpha_i) \in K$ where 

$$\alpha_i = \left\{\begin{array}{cc}0 & \mbox{ if } \rho_i(b) =  0 \\  \rho_i(c)/\rho_i(b) & \mbox{ if } \rho_i(b) \ne 0 \\\end{array}\right. .$$
 
One can check that $\alpha^2 = b, \alpha^3 = c \in A$.  Since $A$ is seminormal in the sense of Traverso, $\alpha = (\rho_i(a))$ for some $a \in A$.  Then, $a^2 = b, a^3 = c$, as desired.
 
 The converse was observed to be true in the preceding paragraph. \qedbox

From this point on, when we speak of a seminormal ring we mean a ring that is seminormal in the sense of Swan.  We wish to define the seminormalization  ${}^+A$ of a reduced ring $A$ to be a seminormal ring   such that $A \subset {}^+A$ is a subintegral extension.  It turns out that any such extension has a universal mapping property that makes it essentially unique; we will make this precise after stating Swan's theorem.
In light of Lemma \ref{swan sn = traverso sn} and Theorem \ref{Traverso char}, if   the ring $A$ is reduced and its total ring of quotients $K$ is a product of fields, then we may define ${}^+A = {}^+_B A$, where $B$ is the normalization of $A$.  If the total quotient ring isn't a product of field, this doesn't always produce a seminormal ring as the following example that appeared in Swan \cite{Swan:1980kn}  illustrates.

\begin{example}  A ring that is seminormal in its total quotient ring need not be seminormal.  Let $(B, \n , k)$ be a local ring and set $A = B[X]/(\n X, X^2) = B[x]$.  Notice that $A$ is local with unique maximal ideal 
$\m = (\n, x)$.  Since $\m$ consists of zero divisors for $A$ the ring $A$ is equal to its own total ring of quotients $K$ and hence $A = {}^+_K A$.  However $A$ is not a seminormal ring, since it isn't reduced.
\end{example}

Swan constructed the seminormalization of a reduced ring using elementary subintegral extensions and mimicking the construction of the algebraic closure of a field.  He also proved the following universal mapping property.
\index{universal mapping property}

\begin{theorem} \cite[Theorem 4.1]{Swan:1980kn}  \label{Swan seminormalization} Let $A$ be any reduced ring.  Then there is a subintegral extension $A \subset B$ with $B$ seminormal.  Any such extension is universal for maps of $A$ to seminormal rings:  If $C$ is seminormal and $\varphi \colon A \to C$, then $\varphi$ has a unique extension $\psi \colon B \to C$.  Furthermore $\psi$ is injective, if $\varphi$ is.
\end{theorem}

Thus if there exist   subintegral extensions $A \subset B$ and $A \subset C$ of a reduced ring $A$ with $B$ and $C$ seminormal, then there is a unique ring homomorphism $\psi \colon B \to C$ which is the identity map on $A$ and $\psi$ is an isomorphism.  Thus in this strong sense, $A$ admits an essentially unique seminormalization.

Mindful of Hamann's examples in \cite{Hamann:1975hq}, Swan calls an extension of rings $A \subset B$ \textit{\textbf{p}}-\textbf{seminormal} if $b \in B, b^2, b^3, pb \in A$ imply $b \in A$; here $p$ denotes  a positive integer.   Swan proves in  \cite[Theorem 9.1]{Swan:1980kn} that a reduced commutative ring is steadfast if and only if it is $p$-seminormal for all rational primes $p$.  This assertion had been proven for domains by Asanuma \cite{MR578866}.  Greither \cite[Theorem 2.3]{MR623811} proved that a projective algebra in one variable over a seminormal ring is a symmetric algebra of an $R$-projective module. 
\index{$p$-seminormal extension}  \index{Swan, Richard}
Before we conclude this section we'd like to mention some results of Greco and Traverso on faithfully flat homomorphisms, pull backs, and completions.   
\index{Greco, Silvio}  \index{Traverso, Carlo}
\begin{defn} A \textbf{Mori ring} is a reduced ring $A$ whose integral closure $\overline{A}$ in its total rings of quotients is a finite $A$-module.
\end{defn}
We point out that an affine ring is Mori.  More generally, a reduced excellent ring is Mori.
\index{Mori ring}   \index{excellent ring}
\begin{theorem} \label{FF result} Consider a  faithfully flat ring homomorphism $f \colon A \to A'$.   Let $B$ be a finite overring of $A$ and put $B' = A' \otimes _ A B$. 
\begin{enumerate} 
\item  If $A'$ is seminormal in $B'$, then $A$ is seminormal in $B$. 
\item If $A'$ is a seminormal Mori ring, then $A$ is a seminormal Mori ring.
\end{enumerate}
\end{theorem}

\noindent \textbf{Proof.}
See \cite[Theorem 1.6 and Cor. 1.7]{MR571055}. \qedbox
\index{faithfully flat homomorphism} \index{pull back} \index{completion}
  
Greco and Traverso also proved this result on pull-backs in order to prove a result for reduced base change.
\index{Greco, Silvio!base change result}  \index{Traverso, Carlo!base change result}
\begin{lemma} \cite[Lemma 4.2]{MR571055} Let $R$ be a ring and let 
\[
\xymatrix{
A \ar[r]^{f} \ar[d]_g & B \ar[d]^u \\
C \ar[r]^v & D}
\]
be a pull-back diagram of $R$-algebras, that is, $$A = B \times_D C := \{ (b,c) \in B \times C \mid u(b) = v(c) \}.$$ Assume that the horizontal arrows are injective and finite, and that the vertical ones are surjective.  If $C$ is seminormal in $D$, then $A$ is seminormal in $B$.
\end{lemma}

This lemma is used to prove the following base change result.  Recall that a homomorphism $f \colon A \to A'$ is \textbf{reduced} if it is flat and all of its fibers $A' \otimes _A \kappa(P)$ are geometrically reduced.  The latter is equivalent to saying that for any prime $P$ of $A$ and any finite extension $L$ of $\kappa(P)$, the ring $B \otimes_A L$ is reduced.                                                                                                                                                                                                                                             
\index{reduced homomorphism}             
             
\begin{theorem}\cite[Theorem 4.1]{MR571055} Let $A \subset B$ be a finite integral extension and let $A \to A'$ be a reduced homomorphism.  If $A$ is seminormal in $B$, then $A'$ is seminormal in $B' = A' \otimes_A B$. 
\end{theorem}

We conclude this section by recalling a result on seminormality and completions.  A proof of this result for algebro-geometric rings can be found in Leahy-Vitulli \cite[Theorem 1.21]{Leahy:1981ta}.
\index{Vitulli, Marie A.} \index{Leahy, John V.}

\begin{theorem}\cite[Cor. 5.3]{MR571055}  An excellent local ring is seminormal if and only if its completion is seminormal.
\end{theorem}

\subsection{Seminormality and the Chinese Remainder Theorem} \label{subsec:CRT}

Let $A$ be a commutative ring and $I_1, \ldots, I_n$ be ideals of $A$.  We say that the  Chinese Remainder Theorem (CRT) holds for $\{ I_1, \ldots , I_n \}$ provided that given $a_1, \ldots, a_n \in A$ such that $a_i \equiv a_j\mbox{ mod }I_i + I_j (i \ne j)$, there exists  an element $a \in A$ such that  $a \equiv a_i \mbox{ mod }  I_i$ for all $i$.
\index{Chinese Remainder Theorem (CRT)}

  Notice that if the ideals $I_j$ are co-maximal we get the statement of  the classical Chinese Remainder Theorem.
  
Note too that in the generalized sense the CRT holds for any pair of ideals since the sequence 
 \[
0 \rightarrow A/(I_1 \cap I_2) \stackrel{\alpha}{\to}  A/I_1 \times A/I_2    \stackrel{\beta}{\to}  A/(I_1 + I_2) \to 0
\]
where $\alpha(a + I_1 \cap I_2) = (a + I_1, a+ I_2)$ and $\beta(a_1+I_1, a_2 + I_2)=a_1 - a_2 + I_1 + I_2$ is always exact.
 
  The generalized Chinese Remainder Theorem need not hold for 3 ideals.  Consider the ideals $I_1 = (X), I_2 = (Y), I_3 = (X-Y) \subset \IC[X,Y]$.   There isn't an element $f \in \IC[X,Y]$ such that $f \equiv Y \mbox{ mod } (X), f \equiv X \mbox{ mod } (Y) \mbox{ and } f \equiv 1 \mbox{ mod }(X - Y)$. 
  
  We now mention a few of the results on seminormality, weak normality,
    and results in the spirit of the Chinese Remainder Theorem that are due to Dayton, Dayton-Roberts and Leahy-Vitulli. \index{Dayton, Barry H.!CRT}\index{Roberts, Leslie G.!CRT}\index{Vitulli, Marie A.!CRT}\index{Leahy, John V.!CRT} These results appeared at about the same time.  The results in Leahy-Vitulli were stated for algebraic varieties over an algebraically closed field of characteristic 0.  The first result appeared both in \cite[Theorem 2]{MR618300} and in \cite[Theorem A]{MR618301}, where a direct proof of the result is given.
 
 \begin{thm}  \label{thm on CRT and sn} Let $A$ be a commutative ring and $I_1, \ldots, I_n$ be ideals of $A$. Suppose that $I_i + I_j$ is radical for $i \ne j$ and let $B = \prod A/I_i$.  Then, the CRT holds for $\{ I_1, \ldots , I_n \}$  if and only if $A/\cap I_i$ is seminormal in $B$.
 \end{thm}
 
 Another form of the CRT appears in \cite[Theorem 2]{MR618300}.  To put this result in some context we remind the reader that if $A \subset B$ is a seminormal extension of rings, then the conductor $A:B$ is a radical ideal of $A$.  A general ring extension $A \subset B$ is seminormal if and only if $A:A[b]$ is a radical ideal of $A[b]$ for every $b \in B$.  Both statements about conductors appeared in Proposition \ref{prop conductors}.  The theorem below is preceded by a lemma on conductors.
 
 \begin{lemma} Let $A$ be a commutative ring and $I_1, \ldots, I_n$ ideals of $A$ with  $\cap I_i = 0$.  Let $B = \prod A/I_i$ and set $\mathfrak{c} = A:B$. Then,
$ \mathfrak{c} = \sum_i \left( \cap_{j \ne i} I_j \right) = \cap_i \left(I_i + \cap_{j\ne i} I_j \right).$

 \end{lemma}
 
 \begin{thm}  \label{big thm on CRT and sn} Let $A$ be a commutative ring and $I_1, \ldots, I_n$ ideals of $A$ with  $\cap I_i = 0$.  Let $B = \prod A/I_i$, $\mathfrak{c} = A:B$ and $\mathfrak{c}_i$ the projection of $\mathfrak{c}$ in $A/I_i$. Assume also that $A/I_i$ is seminormal for each $i$ and that $I_j + I_k$ is radical all $j,k$. Then conditions 1. and 2. below are equivalent.
 \begin{enumerate}
 \item $A$ is seminormal.
 \item The CRT holds for $\{ I_1, \ldots , I_n \}$. 
  \end{enumerate}

\noindent These imply the next 5 conditions, which are equivalent.

  \begin{enumerate}
    \setcounter{enumi}{2}
 \item For each $i$, $\cap_{j \ne i} (I_i + I_j) = I_i + \cap_{j \ne i} I_j$.
 \item For each $i$, $I_i + \cap_{j \ne i} I_j$ is a radical ideal in $A$.
 \item For each $i$, $\mathfrak{c}_i = \cap_{j \ne i} ((I_i + I_j)/I_i)$.
 \item For each $I$, $\mathfrak{c}_i $ is a radical ideal in $A/I_i$.
 \item $\mathfrak{c}$ is a radical ideal in $B$.
 \end{enumerate}
\noindent The above imply the next 2 conditions, which are equivalent. 
\begin{enumerate}
    \setcounter{enumi}{5}
 \item $\mathfrak{c}$ is a radical ideal in $A$.
 \item $\mathfrak{c} = \cup _{j \ne k} (I_j + I_k)$.
 \end{enumerate}
 \noindent If, in addition, $A/\cap_{j \ne i} I_j$ is seminormal for any $i$ then $7. \Rightarrow 1.$
 \end{thm}

There are various results in the spirit of the Chinese Remainder Theorem given in \cite[Section 2]{Leahy:1981ta}, where they are stated for varieties over an algebraically closed field of characteristic 0.  These results appeared at the same time as the algebraic results cited above. Their proofs use an algebro-geometric characterization of weakly normal varieties that will be introduced in Section \ref{subsec: geom props}.

\subsection{Seminormality and the Traverso-Swan Picard Group Result with Coquand's Simplification}\label{subsect:sn and picard grp}
Recall that $\pic (A)$ denotes the Picard group of $A$, i.e., the group of isomorphism classes of rank one projective modules over $A$.  
Given a homomorphism of rings $ \phi \colon A \to B$ we get a homomorphism of groups
\begin{equation}
\pic(\phi) \colon\pic(A) \to \pic (B),
\end{equation}
defined by  extending scalars to $B$.   When $B = A[X]$ and $\phi$ is the natural inclusion, evaluation at zero defines a retract $\rho \colon A[X] \to A$.  We thus have group homomorphisms
\index{Picard group}

\begin{equation}
\pic(A) \stackrel{\pic(\phi)}{\to} \pic (A[X])  \stackrel{\pic(\rho)}{\to} \pic (A),
\end{equation}
whose composition is the identity map. 
We point out that the first homomorphism is always injective and the second always surjective.  The maps $\pic (\phi)$ and $\pic( \rho)$ are isomorphisms if and only if the first is surjective if and only if the second is injective.  Therefore to show that $\pic(\phi)  \colon\pic(A) \to \pic (A[X])$ is an isomorphism  is equivalent to showing: if $P$ is a rank one projective over $A[X]$ and $P(0) :=\pic(\rho)(P) $ is a free $A$-module, then $P$ is a free $A[X]$-module. 

Various people have proven in different cases that $\pic(\phi)$ is an isomorphism if and only if  $A$ is seminormal.  This was first proven by Traverso in  case $A$ is a reduced Noetherian ring with finite normalization \cite[Theorem 3.6]{Traverso:1970zv} and by Swan for arbitrary reduced rings \cite[Theorem 1]{Swan:1980kn}. Both Traverso and Swan used standard K-theoretic results to prove the `if' direction and a construction of Schanuel to prove the `only if' direction.  A recent paper by Coquand \cite{MR2264145} 
 simplified the connection between seminormality and the Picard group result by giving a self-contained proof of Swan's result. For a finite integral extension of reduced rings $A \subset B$, Coquand's work suggests an algorithmic approach to finding a sequence of elements  $a_1, a_2, \ldots , a_n \in B$ such that $a_{i+1}^2, a_{i+1}^3 \in A[a_1, \ldots , a_i] \mbox{ for } i = 1, \ldots , n - 1$ and ${}^+_B A = A[a_1, \ldots a_n]$.  Another algorithm was given by Barhoumi and Lombardi \cite{Barhoumi:2008aa}.
\index{Coquand, Thierry!Picard group} \index{Coquand, Thierry!algorithmic approach to seminormalization}
\index{Barhoumi, Sami!algorithmic approach to seminormalization}  
\index{Lombardi, Henri!algorithmic approach to seminormalization}
\index{Swan, Richard!Picard group}  \index{Traverso, Carlo!Picard group}  
Since $\pic(A) \cong \pic(A_{\mathrm{red}})$ and $A$ seminormal implies $A$ is reduced we may and shall assume that our base ring $A$ is reduced in the remainder of this section.

As pointed out by Swan \cite{Swan:1980kn} one can replace finitely-generated projective modules by ``projection matrices," as we now explain.  If 
$$0 \to K \to A^n  \stackrel{\pi}{\to} P \to 0$$
is a presentation of a finitely-generated projective $A$-module and $\rho \colon P \to A^n$ is a section of $\pi$ then $\rho \circ \pi \colon A^n \to A^n $ is given by an $n \times n$ idempotent matrix $M$ such that $\mathrm{im}(M) = P$; this matrix is referred to as a projection matrix. 

Suppose that $M, M'$ are two idempotent matrices over the ring $A$, not necessarily of the same size. We write $M \cong_A M'$ to express that $M$ and $M'$ present isomorphic $A$-modules.  We write $M \cong_A 1$ to express that $M$ presents a free $A$-module.   Let $P_n$ denote the $n \times n$ matrix [$p_{i j}$] such that $p_{11} = 1$ and all other entries are 0.  Let $I_n$ denote the $n \times n$ identity matrix.

Coquand shows that given  an $n \times n$ idempotent matrix $M$ with entries in $A[X]$, one has  $M(0) \cong_A 1$ only if $M \cong_{A[X]} 1$.  He does this by  establishing a sequence of straightforward lemmas that culminate in the main result.  He first shows that for $A$ seminormal, $M \in A[X]^{n \times n}$ a rank one projection matrix such that $M(0) = P_n$, it must be the case that $M \cong_{A[X]} 1$.  Finally he shows that for   $A$ seminormal, $M \in A^{n \times n}$ a rank one projection matrix such that $M(0) \cong_A 1$, it must be the case that $M \cong_{A[X]} 1$.  The reader is referred to Coquand's paper for the proofs of these results. To give a flavor of the statements of these results and to present Schanuel's example, we mention 2 of Coquand's results.

\begin{lemma} \cite[Lemma 1.1]{MR2264145} 
Let $M$ be a projection matrix of rank one over a ring $A$.  We have $M \cong _A 1$ if and only if there exist
 $x_i, y_j \in A$ such that $m_{i j} = x_iy_j$.  If we write $x$ for the column vector $(x_i)$ and $y$ for the row
  vector $(y_j)$ this can be written as $M = xy$.  Furthermore the column vector $x$ and the row vector $y$ are uniquely defined up to a unit by these conditions: if we have another column vector $x'$ and row vector 
  $y'$ such that $M = x'y'$, then there exists a unit $u \in A$ such that $x = ux'$ and $y' = uy$.
\end{lemma}

Coquand then proves a  result which is well known to many,  that asserts that if $f, g \in A[X]$, $A$ reduced, and $fg = 1$ then $f = f(0), g = g(0)$ in $A[X]$.   Next he proves the following.

\begin{lemma} \cite[Corollary 1.3]{MR2264145}\label{Coquand Cor1.3}
Let $A'$ be an extension of the reduced ring $A$.  Let $M$ be an $n \times n$ projection matrix over $A[X]$ such that $M(0) = P_n$.  Assume that $f_i, g_j \in A'[X]$ are such that $m_{ij} = f_ig_j$ and $f_1(0) = 1$.  If $M \cong_{A[X]} 1$ then $f_i, g_j \in A[X]$.

\end{lemma} 

We state the main theorem of Coquand's paper, which was previously established by Swan. Notice that this result
gives a direct proof for rank one projective modules of the theorem of Quillen-Suslin settling Serre's Conjecture.

\begin{theorem} \cite[Theorem 1]{Swan:1980kn}\cite[Theorem 2.5]{MR2264145} 
If $A$ is seminormal then the canonical map   $\pic(A) \to \pic(A[X_1, \ldots , X_n])$ is an isomorphism.

\end{theorem}

It is straightforward to see, using a construction of Schanuel that originally appeared in a paper by Bass \cite{MR0140542}, more recently appeared in \cite[Appendix A.]{MR2264145}, and in different garb in both \cite{Traverso:1970zv} and \cite{Swan:1980kn}, that if $\pic(\phi) \colon \mathrm{Pic}(A) \to \mathrm{Pic}(A[X])$ is an isomorphism, then $A$ is seminormal.   We give a quick proof.

\begin{lemma}  If $A$ is a reduced ring and $\pic(\phi) \colon \mathrm{Pic}(A) \to \mathrm{Pic}(A[X])$ is an isomorphism, then $A$ is seminormal. 

\end{lemma}

\noindent \textbf{Proof.}

Suppose that $b, c \in A$ and $b^3 = c^2$ and let $B$ be a reduced extension of $A$ containing an element $a$ such that $b = a^2, c = a^3$.  Consider the following polynomials with coefficients in $B$. 

\begin{equation}
f_1 = 1+aX, \; f_2 = bX^2, \; g_1 = (1-aX)(1+bX^2), \; g_2 = bX^2
\end{equation}

The matrix $M = (f_ig_j)$ is a projection matrix of rank one over $A[X]$ such that $M(0) = P_2$.  Assuming that $\pic(\phi) $ is an isomorphism we may conclude that this matrix presents a free module over $A[X]$.  By  Lemma \ref{Coquand Cor1.3}  this implies $f_i, g_j \in A[X]$ and hence $a \in A$. \qedbox

To conclude this section we would like to mention related work of J. Gubeladze \cite{Gubeladze1989:aa} that settled a generalized version of Anderson's Conjecture \cite{AndersonDF:1982}, which  postulated that finitely-generated projectives are free over normal monomial subalgebras of $k[X_1, \ldots, X_n]$. In what follows, by a \textbf{seminormal monoid} we mean a commutative, cancellative, torsion-free monoid $M$ with total quotient group $G$ such that whenever $x \in G$ and  $2x, 3x \in M$ we must have  $x \in M$.  
\index{Gubeladze, Joseph!Anderson's Conjecture}
\index{Anderson's Conjecture}  \index{seminormal monoid}
\begin{theorem}\cite{Gubeladze1989:aa}
Let $R$ be a PID and $M$ be a commutative, cancellative, torsion-free, seminormal monoid.  Then all finitely-generated projective $R[M]$-modules are free. 
\end{theorem}

The original proof of Gubeladze was geometric in nature. Swan gave an algebraic proof of Gubeladze's result in \cite{MR1144038}.  To read more about Serre's Conjecture and related topics the reader is referred to the recent book by Lam \cite{MR2235330}.
\subsection{Seminormal Local Rings in Dimension One} \label{subsect:1 diml sn rings}

A basic fact about a reduced, Noetherian 1-dimensional local ring $A$ is that $A$ is normal if and only if it is a discrete rank one valuation ring.  This leads to a well-known characterization of the normalization of a Noetherian integral domain $A$ with quotient field $K$.  Namely, $b \in K$ is in the integral closure of $A$ in $K$ if and only if $b$ is in every valuation subring $V$ of $K$ containing $A$.

Reduced, Noetherian 1-dimensional seminormal local rings with finite normalization are also fairly well-behaved, particularly in the algebro-geometric setting.  Let $e(R)$ and $\mathrm{emdim}(R)$ denote the multiplicity and embedding dimension of a local ring $(R,\m)$, respectively. The following theorem was proven by E. D. Davis. 
\index{Davis, Edward D.! 1-dimensional seminormal rings}
\begin{theorem} \cite[Theorem 1]{MR0453748}
Let $(R, \m)$ be a reduced, Noetherian, 1-dimensional local ring with finite normalization $S$. The following are equivalent.
\begin{enumerate}
\item $R$ is seminormal.
\item $\mathrm{gr}_{\m}(R)$ is reduced and $e(R) = \mathrm{emdim}(R)$.
\item $\mathrm{Proj}(\mathrm{gr}_{\m}(R))$ is reduced and  $e(R) = \mathrm{emdim}(R)$.
\end{enumerate}
\end{theorem}

If the local ring comes from looking at an algebraic curve at a closed point one can say more.  The following result of Davis generalizes earlier results of Salmon \cite{MR0251032} for plane curves over an algebraically closed field and Bombieri \cite{MR0429874} for arbitrary curves over an algebraically closed field.

\begin{cor}\cite[Corollary 1]{MR0453748} \label{cor bombieri}
Let $x$ be a closed point of an algebraic (or algebroid) curve at which the Zariski tangent space has dimension $n$. Then:
\begin{enumerate}
\item  $x$ is seminormal if, and only if, it is an $n$-fold point at which the projectivized tangent cone is reduced.
\item For an algebraically closed ground field, $x$ is seminormal if, and only if, it is an ordinary $n$-fold point (i.e., a point of multiplicity $n$ with $n$ distinct tangents).
\end{enumerate}
\end{cor}
\index{Bombieri, Enrico!seminormality of ordinary $n$-fold point}
\index{ordinary $n$-fold point}
Thus analytically, a 1-dimensional algebro-geometric seminormal local ring over an algebraically closed field $K$ looks like $K[X_1, \ldots , X_n]/(X_iX_j \mid i \ne j)$.  Here are some specific examples.

\begin{ex}  Let $K$ be an algebraically closed field.  Consider the curve $C = \{ (x,y) \in \A^2 \mid xy - x^6 - y^6 = 0 \}$ given below.
Let $f = xy - x^6 - y^6 $ and $P = (0,0)$.  Since $P$ is a point of multiplicity 2 with 2 distinct tangent lines (the $x$- and $y$- axes), it is an ordinary double point or a node and the local ring of $A = K[x,y]/(f)$ at $P$ is seminormal by Corollary \ref{cor bombieri}.  Since $P$ is the only singular point of $A$ we may conclude that $A$ is seminormal.  This curve is drawn below in Figure  \ref{fig:node}.
\end{ex}
\index{node}
\begin{ex}
Let $K$ be an algebraically closed field.  Consider the curve $C = \{ (x,y) \in \A^2 \mid x^2 = x^4 + y^4  \}$ given below.  Let $f = x^2 - x^4 - y^4 $ and $P = (0,0)$.  Since $P$ is a point of multiplicity 2 with 1 tangent line (the $y$-axis) occurring with multiplicity 2, the local ring of  $A = K[x,y]/(f)$ at $P$ fails to be seminormal by Corollary \ref{cor bombieri}. This curve is called a tacnode and is sketched below in  Figure \ref{fig:tacnode}.
\index{tacnode}
\begin{figure}[h!]
\centering
\subfloat[Node]{\label{fig:node} \includegraphics[width=2in]{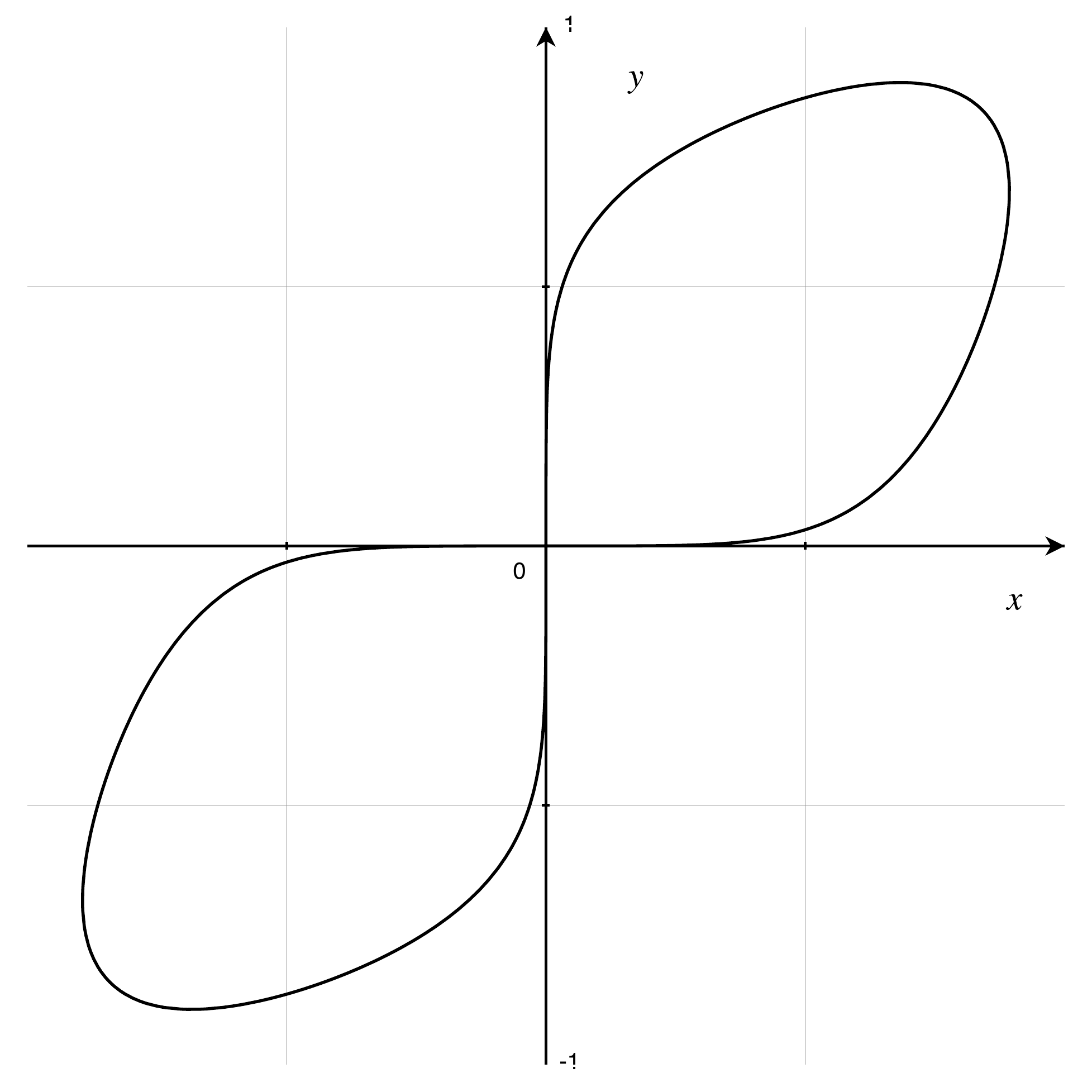}} 
\qquad
\subfloat[Tacnode]{\label{fig:tacnode} \includegraphics[width=2in]{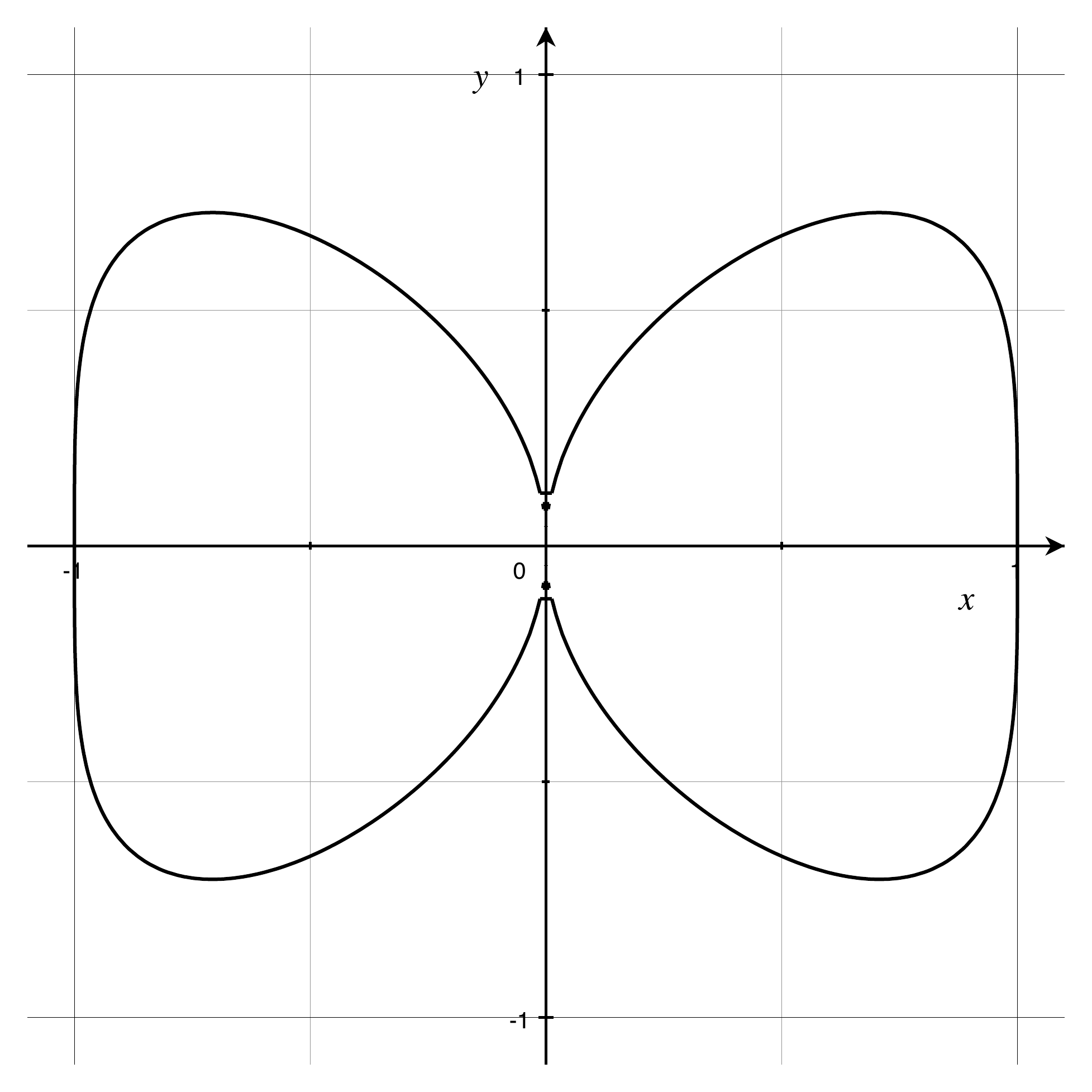}}
\caption{A Seminormal and a nonseminormal curve}
\label{fig:Curves}
\end{figure}
\end{ex}

\section{Weak Normality and Weak Normalization} \label{sect:weak normality and weak normalization}

Let's change gears and speak about weak subintegrality and weak normalization.  
The reader should keep in mind that weak normalization and seminormalization are \emph{not} the same concept in positive characteristic.  

\subsection{ Weakly Subintegral and Weakly Normal Extensions and Weak Normalization Relative to an Extension} \label{subsect:wsi and wn extensions}

We now recall the Andreotti-Bombieri definition of  weak normalization and its first properties.  Recall that the \textbf{characteristic exponent} $e$ of a field $k$ is 1 if $\mathrm{char}(k) = 0$ and is $p$ if $\mathrm{char}(k) = p >0$.  As earlier we only distinguish between a prime ideal $P = P_x$ of a ring $B$ and the corresponding point $x \in \spec(B)$ when it is convenient to do so.
\index{characteristic exponent}
\begin{lemma}  \cite[Prop. 2]{Andreotti:1969aa} \label{lemma on weak local gluing}
Let $(A, \m, k)$ be a local ring,  $e$ be the characteristic exponent of $k$, and $A \subset B$ an integral extension.  Set $$A' = \{ b \in B \mid b^{e^m} \in A + R(B) \mbox{ for some } m \ge 0\}$$ and let $\m' = R(A')$.   For each maximal ideal $\n _i$ of $B$ let $\omega _i \colon A/ \m \to B/ \n_i$ be the canonical homomorphism. The following assertions hold.
\begin{enumerate}
\item  $(A', \m', k')$ is a local ring and the induced extension of residue  fields \\  $k = \kappa(x) \subset k' = \kappa(x')$ is purely inseparable;
\item  $A'$ is the largest intermediate ring that is local and has residue class field a purely inseparable extension of  $k$.
\item An element $b \in B$ is in $A'$ if and only if there exists some integer $m \ge 0$ such that
\begin{enumerate}
\item$ b(x_i)^{e^m} \in \omega_i(\kappa(x)) \mbox{ for all closed points } x_i $ of $\spec(B)$, and
\item $\omega _i^{-1} (b(x_i)^{e^m}) = \omega _j^{-1} (b(x_j)^{e^m}) \mbox{ for all } i, j $.
\end{enumerate}
\end{enumerate}
\end{lemma}

\noindent \textbf{Proof.}
 Either see the proof of \cite[Prop. 2]{Andreotti:1969aa} or modify the proof of Lemma \ref{local gluing} above. \qedbox

\begin{defn} \label{defn weak gluing} We will refer to the ring $A'$ described above as the \textbf{weak gluing} of $A$ in $B$ over $\m$. Letting $x \in \spec (A)$ denote the point corresponding to the maximal ideal $\m$, this ring is sometimes denoted by ${}^*_x A$.
\end{defn}
\index{weak gluing}
\index{${}^*_x A$: weak gluing of $A$ over $x$}
The weak normalization of a ring $A$ in an integral extension $B$ is defined analogously to the seminormalization.

\begin{definition}  \label{defn weak normalization}   Let $A \subset B$ be an integral extension of rings. We define the \textbf{weak normalization} ${}^*_B A$ of $A$ in $B$ to be
\begin{equation}
{}^*_B A = \{ b \in B \mid \forall x \in \spec(A), \exists m \ge 0 \mbox{ such that } (b_x)^{e^m} \in A_x + R(B_x)  \}.
\end{equation}
In this description, $e$ denotes the characteristic exponent of the residue class field $\kappa(x)$.   
\end{definition}
\index{weak normalization}
\index{${}^*_B A $: weak normalization of $A$ in $B$}
Notice that from the definition we have  $A \subset {}^+_B A  \subset {}^*_B A$ for an arbitrary integral extension $A \subset B$.  

Before stating some fundamental properties of weak normalization we introduce some additional terminology.

\begin{definition}  \label{defn weakly subintegral ext} An integral extension of rings $A \subset B$ is said to be \textbf{weakly subintegral} if the associated map $\spec(B) \to \spec(A)$ is a bijection and induces purely inseparable extensions of the residue class fields. An element $b \in B$ is said to be \textbf{weakly subintegral over \emph{B}} provided that $A \subset A[b]$ is a weakly subintegral extension.
\end{definition}
\index{weakly subintegral extension}
Notice that a weakly subintegral extension of fields is a purely inseparable extension.  A subintegral extension must be weakly subintegral (by definition), but the converse need not be true since a purely inseparable extension of fields is weakly subintegral but not subintegral.  The following lemma is the basis   for what Yanagihara calls an elementary weakly subintegral extension.   We include Yanagihara's proof for the convenience of the reader.

\begin{lemma}\cite[Lemma 3]{Yanagihara:1985os}  Let $A \subset A[b]$ be a simple extension of rings, where $p $ is a rational prime and $pb, b^p \in A$.   Then, $A \subset B$ is weakly subintegral.
\end{lemma}

\noindent \textbf{Proof.}

Let $P \subset A$ be a prime ideal.  First suppose that $p \notin P$.  Then, $b = pb/p \in A_P$, which implies $A_P = B_P$ and the residue fields are isomorphic.  Now suppose that $p \in P$, so that $A/P$ is an integral domain of characteristic $p$.  By the Lying Over Theorem, $PB \cap A = P$. Suppose $Q_1, Q_2$ are prime ideals of $B$ lying over $P$.  Consider the extension $A' := A/P \subset B' := B/PB$ and let $Q'_i = Q_i/PB \, (i=1, 2)$.  Then, $Q'_1$ and $Q'_2$ both lie over the zero ideal  in $A'$.  Notice that $f^p \in A'$ for any element $f \in B'$.   Consider $f \in Q'_1$.  Then, $f^p \in Q'_1 \cap A' = 0 \in Q'_2 \Rightarrow f \in Q'_2$.  Hence $Q'_1 \subset Q'_2 $, which implies $Q_1 \subset Q_2$.  Similarly, $Q_2 \subset Q_1$ and, hence, $Q_1 = Q_2$.  Thus there is a unique prime ideal $Q$ of $B$ lying over $P$ and $\kappa(P) \subset \kappa(Q)$ is purely inseparable. \qedbox

\begin{definition} An \textbf{elementary weakly subintegral extension} is a simple extension $A \subset A[b]$ of rings such that $pb, b^p \in A$, for some rational prime $p$.

\end{definition}
\index{weakly subintegral extension!elementary weakly subintegral extension}
\medskip

The notion of a weakly normal extension, which we now define, is \textit{complementary} to that of a weakly subintegral extension.

\begin{definition} If $A \subset B$ is an integral extension of rings, we say $A$ is \textbf{weakly normal }in $B$ if there is no subextension $A \subset C \subset B$ with $C \ne A$ and $A \subset C$ weakly subintegral.
\end{definition}
\index{weakly normal extension}
The following result is well known and can be proven by modifying the proof of Theorem \ref{Traverso char}. 

\begin{theorem} \label{Yangihara char} Let $A \subset B$ be an integral extension of rings. The following assertions hold.
\begin{enumerate}
 \item The extension $A \subset {}^*_B A$ is  weakly subintegral;
 \item If $A \subset C \subset B$ and $A \subset C$ is weakly subintegral, then $C \subset {}^*_B A$; 
 \item  The extension ${}^*_B A \subset B$ is weakly normal; and
 \item  ${}^*_B A$ has no proper subrings containing $A$ and weakly normal in $B$.
\end{enumerate}
\end{theorem}

The reader will note that ${}^*_BA$ is the unique largest weakly subintegral extension of $A$ in $B$ and is minimal among the intermediate rings $C$ such that $C \subset B$ is weakly normal.  With this result in hand, we see that an integral extension $A \subset B$ is weakly normal if and only if $A = {}^*_B A$.   In light of the inclusions $A \subset {}^+_B A  \subset {}^*_B A$, a weakly normal extension is necessarily seminormal, but the converse is not true as Example \ref{sn ext not wn} will illustrate.

H. Yanagihara proved a weakly subintegral analogue of Swan's characterization of subintegral extensions in terms of extensions of maps into fields, which we now state.

\begin{lemma} \cite[Lemma 1]{Yanagihara:1985os} \label{field char of weakly subintegral} Let $A \subset B$ be an integral extension of rings.  It is a weakly subintegral extension if and only if for all fields $F$ and homomorphisms $\varphi \colon A \to F$ there exists at most one homomorphism $\psi \colon B \to F$ extending $\varphi$.
\end{lemma}

 The positive integer $p$ that appears in the work of Hamann and Swan also appears in  Yanagihara's characterization of a weakly normal extension.
Yanagihara showed \cite{MR714467} that if a ring  $A$ contains a field of positive characteristic, an integral extension $A \subset B$ is weakly normal provided that $A \subset B$ is seminormal and $A$ contains every element $b \in B$ such that $b^p$ and $pb \in A$, for some rational prime $p$. This result has been widely cited but didn't appear in print until much later.   Subsequently,
Itoh \cite[Proposition 1]{MR723066} proved the above-mentioned characterization holds for any integral extension $A \subset B$ and Yanagihara gave an alternate proof \cite[Theorem 1]{MR714467}.   
We now recall the generalized result.

\begin{prop} Let $A \subset B$ be an integral extension of rings.  Then, $A$ is weakly normal in $B$ if and only if the following hold:
\begin{enumerate}
\item  $A$ is seminormal in $B$; and
\item  $A$ contains each element $b \in B$ such that $b^p, pb \in A$, for some rational prime $p$.
\end{enumerate}
\end{prop}

Thus  $A$ is weakly normal in $B$ if and only if $A$ doesn't admit any proper (elementary) subintegral or weakly subintegral extensions in $B$.  In fact, this characterizes a weakly normal extension by \cite[Corollary to Lemma 4]{Yanagihara:1985os}.  

For an integral extension $A \subset B$ the weak normalization ${}^* _ B A$ is the filtered union of all subrings of $B$ that can be obtained from $A$ by a finite sequence of elementary weakly subintegral extensions.  We now offer an example of a seminormal extension of rings that isn't weakly normal.

\begin{ex}\label{sn ext not wn}  Let $K$ be a field of characteristic 2, $X$ be an indeterminate, and consider the  integral extension $A := K[X^2] \subset B := K[X]$.  Since $X^2$ and $2X = 0$ are both in $A$ but $X$ is not in $A$, this extension is not weakly normal.  However, it is a seminormal extension as seen in Example \ref{sn ext ex}.

\end{ex}

Inspired by Swan's intrinsic definition of a seminormal ring, Yanagihara made the following definition in \cite{Yanagihara:1985os}.
\index{Yanagihara, Hiroshi!weakly normal ring}

\begin{definition} \label{Yanagihara intrinsic def wn} A ring $A$ is said to be \textbf{weakly normal} if $A$ is reduced and the following conditions hold:
\begin{enumerate}
\item  For any elements $b, c \in A$ with $b^3 = c^2$, there is an element $a \in A$ with $a^2 = b$ and $a^3 = c$; and
\item  For any element $b, c, e \in A$ and any non-zero divisor $d \in A$ with $c^p = bd^p$ and $pc = de$ for some rational prime $p$, there is an element $a \in A$ with $b = a^p$ and $e = pa$.
\end{enumerate}
\end{definition}
\index{weakly normal ring}
As the reader may recall, condition 1. implies that the ring $A$ is reduced so we may and shall omit that requirement in the definition of a weakly normal ring.  Since $A$ is reduced, the element $a$ that occurs in condition 1. is necessarily unique.   Since this definition is quite complicated we will look at some special instances of weakly normal rings.

\begin{lemma} A product of fields is weakly normal.
\end{lemma}

\noindent \textbf{Proof.} We first show that a field $K$ is weakly normal.  Suppose $0 \ne b, c \in K$ and $b^3 = c^2$.  Setting $a = c/b$ we have $b = a^2, c = a^3$.  Now suppose that $b, c, e \in K$ and  $ 0 \ne d \in K$ with $c^p = bd^p$ and $pc = de$ for some rational prime $p$.   Setting $a = c/d$ we have $b = a^p, \, e = pa$.   Now let $K = \prod K_i$ be a product of fields.  Suppose $b= (b_i), c = (c_i) \in K$ and $b^3 = c^2$.  Setting $a_i = c_i/b_i$ whenever $b_i \ne 0$ and $a_i = 0$ whenever $b_i = 0$, we have $b = a^2, c = a^3$.  Now suppose that $b, c, e \in K$ and  $  d \in K$ a non-zero divisor with $c^p = bd^p$ and $pc = de$ for some rational prime $p$.  Since $d = (d_i)$, where each $d_i \ne 0$ we can let $a_i = c_i/d_i$ and set $a = (a_i)$.  Then, $b = a^p, \, e = pa$. \qedbox

\begin{lemma} \label{normal is wn}

Let A be a reduced ring whose total quotient ring is a product of fields.  If $A$ is normal, then $A$ is weakly normal.

\end{lemma}

\noindent \textbf{Proof.}
Let $A$ be a normal ring whose total quotient field is a product of fields.  Since a ring $A$ with such a total quotient ring is seminormal if and only if it is equal to ${}^+_B A$, where $B$ is the normalization of $A$, we may conclude that $A$ is seminormal.  Now suppose that $b, c, e \in A$ and  $ d \in A$ is not a zero divisor such that $c^p = bd^p$ and $pc = de$ for some rational prime $p$.  Since $(c/d)^p = b \in A$ and $A$ is normal, we must have $c/d \in A$ so we may let $a = c/d$ as in the above example. \qedbox
 
 We now recall a pair of results by Yanagihara that we will find helpful.
 
 \begin{prop}\cite[Propositions 3 and 4]{Yanagihara:1985os} \label{Yanag props} Let $A$ and $B$ be reduced rings.
\enumerate
\item  If $A$ is weakly normal and $B$ is a subring of the  total quotient ring of $A$ containing $A$, then $A$ is weakly normal in $B$.  
 \item  If $A$ is a subring of  a weakly normal ring $B$ such that any non-zero divisor in $A$ is also not a zero divisor in $B$ and $A$ is weakly normal in $B$, then $A$ is weakly normal.
 \end{prop}
 
 The next result is a weakly normal analog of the fact that a ring whose total quotient ring is a product of fields  is a  seminormal ring if and only if it is equal to its seminormalization in $\overline{A}$.   
 
 \begin{cor}  \label{wn = wn in normalization} Let $A$ be a reduced ring whose total quotient ring is a product of fields and let $B$ denote the normalization of $A$.  The following are equivalent.
 \begin{enumerate}
 
\item $A$ is weakly normal;
\item $A = {}^*_B A$; and 
\item $A$ is weakly normal in $B$.
 \end{enumerate}
 \end{cor}
 
\noindent \textbf{Proof.}
The last two items are equivalent as remarked after the proof of Theorem \ref{Yangihara char}.  The equivalence of the first and last statements follows from Lemma \ref{normal is wn} and  the second part of \ref{Yanag props}. \qedbox
 
When the base ring is reduced and contains a field of positive characteristic,  weak normality has a simple characterization.  Notice that to say a reduced ring $A$ is $n$-closed in its normalization for some positive integer $n$ is the same thing as saying that $A$ contains each element $b$ of its total quotient ring such that $b^n \in A$.  If the latter condition holds we will say $A$ is $n$-closed in its total quotient ring $K$ even though $A \subset K$ need not be an integral extension.

 \begin{cor} Suppose that $A$ is a reduced ring whose total quotient ring is a product of fields and that $A$ contains a field  of characteristic $p > 0$.  Then, $A$ is weakly normal if and only if $A$ is $p$-closed in its total quotient ring.
 \end{cor}
 
\noindent \textbf{Proof.}

 Let $A$ be a reduced ring whose total quotient ring $K$ is a product of fields and suppose that $A$ contains a field $k$ of characteristic $p > 0$.   By Lemma \ref{wn = wn in normalization}  and Hamman's criterion it suffices to show that $A$ is (2,3)-closed in $K$ and contains each element $b \in K$ such that $qb, b^q \in A$ for some rational prime $q$.  Suppose that $b \in K$ and $qb, b^q \in A$ for some rational prime $q$.  If $q \ne p$, then $b = qb/q \in A$.  Suppose $q = p$.  Since $A$ is $p$-closed in $K$ by assumption, $b \in A$.  Now suppose that $b \in K$ and $b^2, b^3 \in A$.  Then, $b^n \in A$ for all $n \ge 2$.  In particular, $b^p \in A$ and hence $b \in A$. \qedbox
 
 Looking at the complementary notion of weak subintegrality in positive characteristic we have another simple characterization, which we now recall.
  
 \begin{prop}\cite[Theorems 4.3 and 6.8] {Reid:1996dw} Let $A \subset B$ be an extension of commutative rings and suppose that $A$ contains a field of characteristic $p > 0$. Then $b \in B$ is weakly subintegral over $A$ if and only if $b^{p^n} \in A$ for some $n \ge 0$. 
 
 \end{prop} 
 
 We mention a few more results of Yanagihara regarding localization, faithfully flat descent, and pull backs.
  \index{Yanagihara, Hiroshi!faithfully flat descent and weak normalization}
    \index{Yanagihara, Hiroshi!pull backs and weak normalization}  \index{faithfully flat homomorphism} \index{pull back}
    
\begin{prop}\cite[Prop. 7]{Yanagihara:1985os}
A ring $A$ is weakly normal if  all localizations of $A$ at maximal ideals are weakly normal.  
\end{prop}

\begin{prop}\cite[Prop. 5]{Yanagihara:1985os}
Let $B$ be a weakly normal ring and $A$ be a subring of $B$ such that $A$ has only a finite number of minimal prime ideals. If $B$ is faithfully flat over $A$, then $A$ is weakly normal.
\end{prop}

\begin{prop} \cite[Prop. 1]{MR577861} Let 
\[
\xymatrix{
D \ar[r]^{\alpha} \ar[d]_\beta & A \ar[d]^f \\
B \ar[r]^g & C}
\]
be a pull-back diagram of commutative rings.  Assume that $A$ is weakly normal domain, that $B$ is a perfect field of positive characteristic $p$, and that $C$ is reduced.  Then, $D$ is weakly normal.
\end{prop}

Yanagihara does not construct the weak normalization of a general commutative ring.

Itoh went on to prove a result like Greither's.  
He showed \cite[Theorem]{MR723066} that a finite extension $A \subset B$ of reduced Noetherian rings is a weakly normal extension if and only if for every $A$-algebra $C$, $B \otimes_A C \cong \mathrm{Sym}_B(Q), Q \in \pic (B)$, implies $C \cong \mathrm{Sym}_A(P),$ for some $P \in \pic (A)$.

\subsection{Systems of (Weak) Subintegrality} \label{subsect:SOSI}

The reader may have noticed that the only way to determine if
 an element $b$ in an integral extension ring $B$ of a ring $A$ is weakly subintegral over $A$ is to check if the \emph{ring extension}  $A \subset A[b]$ is weakly subintegral. Beginning in 1993, in a series of papers \cite{Roberts:1993lo}, \cite{Reid:1993ij}, \cite{Roberts:1994qa}, \cite{Reid:1995fy}, \cite{Reid:1995pr}, \cite{Reid:1996dw}, \cite{Reid:1996ys} by Reid, Roberts, and Singh,  a genuine elementwise criterion for weak subintegrality was introduced and developed.  At first they worked with $\Q$-algebras, where weak subintegrality and subintegrality coincide, and discussed a criterion for subintegrality that involved what they called a \textbf{system of subintegrality} or SOSI.  Over time, they were able to handle general rings, but only if they discussed weak subintegrality rather than subintegrality.  However, since the term system of subintegrality had already been used extensively in their earlier papers they kept it instead of redefining the system as a system of \emph{weak} subintegrality.  
 \index{Roberts, Leslie G.!SOSI}
 \index{Reid, Les!SOSI}
 \index{Singh, Balwant!SOSI}
Reid, Roberts and Singh constructed in \cite[section 2]{Reid:1996dw} a ``universal weakly subintegral extension" based on systems of subintegrality.   Later Roberts \cite{MR1792146} wrote a paper that helped to elucidate the earlier work.  Much more recently Gaffney and the current author \cite{GaffVit2008:aa} gave a more intuitive geometric description of the elementwise criterion for weak subintegrality for rings that arise in the study of classical algebraic varieties or complex analytic varieties.  We offer a new algebraic criterion for weak subintegrality in this section.
\index{universal weakly subintegral extension}

\begin{theorem}  \cite[Theorems 2.1 ,5.5, and 6.10]{Reid:1996dw} \label{RRS equiv wn thm}  Let $A \subset B$ be an extension of rings, $b \in B$,  and $q$  a nonnegative integer.  The following statements are equivalent.
\begin{enumerate} 
\item There exists a positive integer $N$ and elements $c_1, \ldots , c_q \in B$ such that 
\[
 b^n + \sum_{i=1}^q {n \choose i} c_ib^{n-i} \in A  \mbox{ for all } n \ge N.
\]
 \item There exists a positive integer $N$ and elements $c_1, \ldots , c_q \in B$ such that 
\[
b^n + \sum_{i=1}^q {n \choose i} c_ib^{n-i} \in A  \mbox{ for } N \le n \le 2N + 2q-1.
\]
  \item There exist elements $c_1, \ldots , c_q \in B$ such that 
\[
 b^n + \sum_{i=1}^q {n \choose i} c_ib^{n-i} \in A  \mbox{ for all } n \ge 1.
\]
\item There exist elements $a_1, a_2, a_3, \ldots \in A$ such that
\[
b^n + \sum_{i=1}^n (-1)^i{n \choose i} a_ib^{n-i} = 0 \mbox{ for all }  n  > q.
\]
\item There exist elements $a_1, \ldots, a_{2q+1}$ in $A$ such that \\
\begin{equation}\label{wsi eqns}
b^n + \sum_{i=1}^n (-1)^i{n \choose i} a_ib^{n-i} = 0 \mbox{ for } q+1 \le n  \le 2q + 1.
\end{equation}
  \item  The extension $A \subset A[b]$ is weakly subintegral.
\end{enumerate}
\end{theorem}

Notice condition 2. in the above theorem is a finite version of condition 1. and that condition 3. is a special case of condition 1. with $N = 1$.  Similarly, condition 5. is a finite version of condition 4.

We now recall the definition that was the basis of the work of Reid, Roberts, and Singh.

\begin{definition}  Let $A \subset B$ be an extension of rings and $b \in B$. A \textbf{system of subintegrality} (SOSI) for $b$ over $A$ consists of a nonnegative integer $q$, a positive integer $N$, and elements $c_1, \ldots , c_q \in B$ such that 
\begin{equation}\label{SOSI defn}
  b^n + \sum_{i=1}^q {n \choose i} c_ib^{n-i} \in A  \mbox{ for all } n \ge N.
\end{equation}
\end{definition}
\index{system of subintegrality (SOSI)}

Originally, an element $b$ admitting a SOSI was called \textbf{quasisubintegral} over $A$ but Reid, Roberts, and Singh later dropped this terminology for reasons which will soon be clear.   We next summarize several results of Reid,  Roberts, and Singh.
\index{quasisubintegral}

In light of the equivalence of conditions 1. and 6. we say an element $b \in B$ is  weakly subintegral over $A$ provided that $b$ admits a SOSI over $A$; we will not mention a quasisubintegral element again.   By this theorem,  an element $b$ is weakly subintegral over $A$ if and only if $b$ satisfies a highly structured sequence (\ref{wsi eqns}) of equations of integral dependence.  



Let $F(T) = T^n + \sum_{i=1}^n (-1)^n{n \choose i} a_iT^{n-i}$.  Notice that $F'(T) = nF_{n-1}(T)$.   So the equations that appear in part 5. of Theorem \ref{RRS equiv wn thm} are rational multiples of the derivatives of  the equation of highest degree.   In the section to follow we will show directly  that this system of equations implies that $b$ is weakly subintegral over $A$.  Hopefully our proof will shed more light on why this system of equations of integral dependence implies weak subintegrality.

\subsection{A New Criterion for Weak Subintegrality} \label{subsect:fresh approach}

We will now develop a new criterion for an element to be weakly subintegral over a subring.   First we look at the case of an extension of fields.

\begin{lemma} \label{wn el over field} Suppose that $K \subset L$ is an extension of fields and that $x \in L$. Then, $x$ is weakly subintegral over $K$ if and only if  $x$ is a root of some monic polynomial $F(T) \in K[T]$ of degree $n > 0$ and its first $\lfloor \frac{n}{2} \rfloor$ derivatives.    
\end{lemma}

\noindent \textbf{Proof.} 
Throughout the proof we let $\ell = \lfloor \frac{n}{2} \rfloor$.  First assume that  $x$ is a root of some monic polynomial $F(T) \in K[T]$ of degree $n > 0$ and its first $\ell$ derivatives.    
Let $f(T)$ denote the minimal polynomial of $x$ over $K$. Then, $F(T) = f(T)g(T)$ for some monic polynomial $g(T) \in K[T]$. Looking at the Taylor series expansion of $F(T)$ in $T-x$ we see that  our hypotheses imply that
\begin{eqnarray*}
F(T) &=& \sum_{k=0}^n \frac{F^{(k)}(x)}{k!} (T-x)^k \\
 &=&  (T-x)^{\ell + 1}G(T) \\
 &=&  f(T)g(T).
\end{eqnarray*}
 Notice that the coefficients of $\frac{F^{(k)}(x)}{k!}$ are actually integer multiples of elements of $K$ and $G(T)$ is monic with coefficients in $K(x)$.

We break the rest of the argument into 2 cases depending on whether the characteristic of $K$ is 0 or positive.

Suppose first that char($K$)=0.  We claim that $\deg(f) = 1$, i.e., $x \in K$.  Suppose to the contrary that $\deg(f) > 1$.     Since $x$ is a simple root of $f(T)$ we may conclude $(T-x)^{\ell} \mid g(T)$.  This implies $g = f^{\ell}g_{\ell}$ for some monic $g_{\ell} \in K[T]$ and, in turn, $F = f^{\ell + 1}g_{\ell}$, which is absurd since the right hand side has degree at least $2(\ell + 1) > n = \deg(F)$. Hence, $\deg(f) = 1$ and $x \in K$ as asserted.

Now assume that char($K) = p > 0$.  Write $f(T) = h(T^{p^m})$ for some polynomial $h(T) \in K[T]$ having only simple roots.  If $m = 0$ proceed as above to deduce that $x \in K$. So assume $m >0$. Let $r = \deg(h)$ so that 
\begin{eqnarray*}
f(T) &=& h(T^{p^m}) \\
&=& (T^{p^m} - z_1) \cdots (T^{p^m} -z_r) \\
&=& (T - x_1)^{p^m} \cdots (T - x_r)^{p^m},
\end{eqnarray*}
over some splitting field for $f$ over $K$ containing $x = x_1$.
We claim that $r = \deg(h) = 1.$    Since $f'(T) = 0$ we must have $F^{(k)}(T) = f(T)g^{(k)}(T)$ for all $k \ge 1$ and hence each $x_i$  is a root of $F(T)$ and its first $\ell$ derivatives, which is absurd since $F$ has at most one root of multiplicity at least $\ell + 1$ by degree considerations. Thus $r = 1$, which implies $x^{p^m} \in K$.  Thus $K \subset K[x]$ is purely inseparable.

Now assume that $x$ is weakly subintegral over $K$.  If $K$ has characteristic 0, then $x \in K$ and $F(T) = T-x$ is the required polynomial.  If $K$ has characteristic $p> 0$ then $x^{p^m} \in K$ for some $m > 0$ and $F(T) = T^{p^m} - x^{p^m}$ is the required polynomial.  \qedbox

An immediate consequence of the above proof is the following.

\begin{cor}  Let $K \subset L$ be an extension of fields and $x \in L$ be algebraic over $K$.  Then, $x$ is weakly subintegral over $K$ if and only if the minimal polynomial of $x$ over $K$ is $f(T) = T^{e^m} - x^{e^m}$ for some positive integer $m$, where $e$ is the characteristic exponent of $K$. 
\end{cor}

Next we present our new element-wise criterion for weak subintegrality and give an elementary self-contained proof that an element satisfying a system of equations  (\ref{wsi eqns})  is weakly subintegral over the base ring of an arbitrary integral extension.   

\begin{prop}  \label{nec for wn} Let $A \subset B$ be an integral extension of rings and  $b \in B$. Then, $b$ is weakly subintegral over $A$ if and only if there is a monic polynomial $F(T)$ of degree $n> 0$ such that $b$ is a root of $F(T)$ and its first $\lfloor \frac{n}{2} \rfloor$ derivatives. 
\end{prop}

\noindent \textbf{Proof.}

First assume that $F(T)$ is a monic polynomial of degree $n> 0$ such that $b$ is a root of $F(T)$ and its first 
$\lfloor \frac{n}{2} \rfloor$ derivatives. 
Replacing $B$ by $A[b]$ we may and shall assume that $A \subset B$ is a finite integral extension.   Let $P \in \spec(A)$ and $e$ be the characteristic exponent of $\kappa(P)$. We will show that $b_P^{e^m} \in A_P + R(B_P)$ for some $m \ge 1$.  Replacing $A, B, \mbox{ and } b$ by $A_P, B_P, \mbox{ and } b_P$, respectively, we may and shall assume that $(A,P,K)$ is local and $K$ has characteristic exponent $e$. Let $Q_i$ denote the maximal ideals of $B$, and $ \omega_i \colon K =\kappa(P) \to \kappa(Q_i)$ the canonical injections of residue class fields.   Let $\overline{F}$ denote the polynomial obtained from $F$ by reducing coefficients modulo $P$.  Fix a maximal ideal $Q_i$ of $B$ and let $b(Q_i)$ denote the image of $b$ in $\kappa(Q_i) = B/Q_i.$ By Lemma \ref{wn el over field} we may conclude $b(Q_i)^{e^m} \in \omega_i(\kappa(P))$ for some $m$.  Since there are finitely many maximal ideals in $B$ there is some positive integer $m$ such  $b(Q_i)^{e^m} \in \omega_i(\kappa(P))$ for all $i$.  Since $\overline{F}$ has at most one root of multiplicity at least $\lfloor \frac{n}{2} \rfloor$ we may conclude that  $\omega _i^{-1} (b(Q_i)^{e^m}) = \omega _j^{-1} (b(Q_j)^{e^m}) \mbox{ for all } i, j$.  Hence $b^{e^m} \in A + R(B)$, as desired. 

Now assume that $b$ is weakly subintegral over $A$.  Then, there exist  $q \ge 0$ in $\Z$ and elements $a_1, \ldots, a_{2q+1}$ in $A$ such that 
$b^m + \sum_{i=1}^m (-1)^i{n \choose i} a_ib^{m-i} = 0$ for all integers $q$ with $q+1 < m  < 2 q+1$ by part 5. of  Theorem \ref{RRS equiv wn thm}.
Let $n = 2q+1$ and $F(T) = T^n + \sum_{i=1}^n (-1)^i{n \choose i} a_iT^{n-i} $.  By our earlier remarks, $b$ is a root of  $F(T)$ and its first $\lfloor \frac{n}{2} \rfloor$ derivatives. \qedbox

We feel that there should be a direct proof of the `only if' part of the above result, i.e., a proof that doesn't depend on Theorem \ref{RRS equiv wn thm} or  \cite[Theorem 6.8]{Reid:1996dw}.

\subsection{First Geometric Properties of Weakly Normal Varieties} \label{subsec: geom props}

In the introduction we traced the history of weak normality and weak normalization, beginning with its roots in complex analytic space theory. We'd like to further discuss the history of weakly normal complex analytic spaces and then take a brief glimpse at the theory of weakly normal algebraic varieties from a geometric point of view.

Weak normalization first was defined for complex analytic spaces by Andreotti and Norguet \cite{Andreotti:1967aa}. Some say that their quotient space construction was already implicit in the work of Cartan \cite{MR0139769}.
The sheaf $\mathcal{O}^c_X$ of $c$-regular functions on a complex analytic space $X$ is defined by setting its sections on an open subset $U \subset X$ to be those continuous complex-valued functions on $U$ which are holomorphic at the regular points of $U$; this coincides with setting $\Gamma(U, \mathcal{O}^c_X)$ equal to the set of continuous complex-valued functions on $X$ that become regular when lifted to the normalization of $X$.  The weak normalization $X^w$ of $X$ represents the sheaf $\mathcal{O}^c_X$ in the the category of complex analytic spaces, i.e., $\mathcal{O}^c_X$  becomes the sheaf of germs of holomorphic functions on $X^w$.
A complex analytic space is weakly normal if and only if the sheaf of $c$-regular functions on $X$ coincides with the sheaf of regular functions on $X$.   
\index{$c$-regular function}  \index{weakly normal complex analytic space}
\index{Andreotti, Aldo!weakly normal complex analytic space}  \index{Norguet, Fran\c{c}ois!weakly normal complex analytic space}
The weak normalization of a what is today called a scheme was introduced by A. Andreotti and E. Bombieri \cite{Andreotti:1969aa}.  After introducing the notion of ``gluing" the prime ideals of $B$ lying over the unique maximal ideal of $A$,
  where $A \subset B$ is an integral extension,  Andreotti and Bombieri constructed the weak normalization of the structure sheaf of a scheme  pointwise using the gluing they previously defined for
   local rings. 
They then turned their attention to defining and constructing the weak normalization $ \sigma \colon X^* \to X$ of a reduced algebraic scheme $X$ over an arbitrary field $K$.      
They showed that $(X^*, \sigma)$ is maximal among all pairs $(Z, g)$ consisting of an algebraic scheme $Z$ over $K$ and a $K$-morphism $g \colon Z \to X$ that is birational and a universal homeomorphism,
   where the latter means that all maps $Z' \to X'$ obtained by base change are homeomorphisms.   
\index{Andreotti, Aldo!weak normalization of scheme} \index{Bombieri, Enrico!weak normalization of scheme}   \index{universal homeomorphism}
If we study algebraic varieties rather than schemes we can either proceed as in the case of schemes or we can capture some of the ideas that permeated the original complex analytic theory of weakly normal spaces.  We will build on the complex analytic viewpoint.

             In this section,  let $K$ be a fixed algebraically closed field of characteristic 0.  When we speak of an algebraic variety over $K$ we assume that the underlying topological space is the set of closed points of a reduced, separated scheme of finite type over $K$. By an affine ring (over $K$) we mean the coordinate ring of an affine variety (over $K$).   By assuming char $K = 0$ we avoid all inseparability problems and hence the operations of seminormalization and weak normalization coincide.   We will use the latter terminology. 
                          
       Consider an algebraic variety $X$ defined over $K$.  One might expect that if you define the sheaf $\mathcal{O}^c_X$ of $c$-regular functions on $X$ so that its sections on an open subset $U \subset X$ are the continuous $K$-valued functions that are regular at the nonsingular points of $U$ then  $\mathcal{O}^c_X$ becomes the sheaf of regular functions on 
the weak normalization of $X$ (see \cite[Definition 2.4]{Leahy:1981ta} and the corrected definition by the current  author \cite[Definition 3.4]{Vitulli:1987ai}).  This isn't the case due to the very special nature of the Zariski topology in dimension one \cite[Example 3.3]{Vitulli:1987ai}.  We recall the corrected definition now.
\index{Vitulli, Marie A.!$c$-regular function}

\begin{defn} \cite[Definition 3.4]{Vitulli:1987ai}  Let $X$ be a variety and let  $\pi: \tilde{X}\to X$ be the normalization of $X$. The sheaf of $c$-regular functions on  $X$, denoted by $\mathcal{O}^c_X$, is defined as follows.  For an open subset $U \subset X$, we let $\Gamma(U,\mathcal{O}^c_X) = \{ \phi \colon U \to K \mid \phi \circ \pi \in \Gamma(\pi^{-1}(U),\mathcal{O}_{\tilde{X}}) \}$.
\end{defn}
Notice that  a $c$-regular function on an affine variety $X$ may be identified with a regular function on the normalization $\tilde{X}$ of $X$ that is constant on the fibers of $\pi \colon \tilde{X} \to X$. This observation and Theorem \ref{c regular}  below show that the $c$-regular functions on an affine variety may be identified with the regular functions on the weak normalization of $X$. We now describe a $c$-regular function without reference to the normalization.   In this result the underlying topology is the Zariski topology.
\index{$c$-regular function}
\begin{thm}  \cite[Theorem 3.9]{Vitulli:1987ai} Let $X$ be an affine  variety without any 1-dimensional components and consider a function $\phi \colon X \to K$.   Then, $\phi$  is $c$-regular if and only if every polynomial in $\phi$ with coefficients in $\Gamma(X, \mathcal{O}_X)$ is continuous and the graph of $\phi$ is closed in $X \times K$.
\end{thm}

We wish to present some characterizations of the weak normalization of an affine ring.  First we recall  a variant of a well-known result.

\begin{lemma}\cite[Theorem 7, p. 116]{MR0366917} Let $\pi \colon Y \to X$ be a dominating finite morphism of irreducible varieties and let $n = [K(Y):K(X)]$ denote the degree of the extension of fields of rational functions.  Then there is a nonempty open subset $U$ of $X$ such that for each $x \in U$ the fiber $\pi^{-1}(x)$ consists of $n$ distinct points.           
\end{lemma}
\index{finite morphism}  \index{dominating morphism}

With this result in hand  one can prove the following result, which has an immediate corollary.

\begin{theorem}\cite[Theorem 2.2]{Leahy:1981ta} \label{c regular} Let $A \subset B$ be a finite integral extension of affine rings and define $A'$ by  $A' = \{b \in B \mid b_x \in A_x + R(B_x) \; \forall x \in X = \mathrm{Var}(A) \}$.  Then $A' = {}^+_BA$.  Thus if $ \pi \colon Y =\mathrm{ Var}(B) \to X$ is the induced morphism, then ${}^*_B A$ consists of all regular functions $f$ on $Y$ such that $f(y_1) = f(y_2)$ whenever $\pi(y_1) = \pi(y_2)$.
\end{theorem}

\index{regular function}
If you work with an integral extension $A \subset B$ of affine  rings over an algebraically closed field of characteristic 0, the condition that the residue fields be isomorphic in the definition of a seminormal extension is redundant, as we shall now see.  Furthermore, you only need to verify that for each maximal ideal of $A$ there exists a unique prime ideal of $B$ lying over $A$.

\begin{cor} Let $A \subset B$ be an integral extension of affine rings over $K$.  The extension is weakly subintegral if and only if the induced map of affine varieties $\mathrm{Var}(B) \to \mathrm{Var}(A)$ is a bijection.
\end{cor}

\noindent \textbf{Proof.}
The only if direction follows from the definition.  Now suppose that the induced map of varieties is a bijection.  Let $Q \in \spec (B)$ and apply the preceding Lemma to the induced map $\mathrm{Var}(B/Q) \to \mathrm{Var}(A/(Q \cap A))$ to deduce that the induced map of fields of rational functions is an isomorphism.  Since $Q$ was arbitrary, the extension is weakly subintegral. \qedbox

\medskip
We now present an example of a weakly normal surface.

\begin{ex}[The Whitney Umbrella]  Let $A = K[u,uv, v^2] \subset B = K[u,v]$, where $K$ is an algebraically closed  field of characteristic 0.  Notice that $B$ is the normalization of $A$.  We claim that $A$ is weakly normal.  Let $\pi: \A^2 \to \A^3$ be given by $\pi(u,v) = (u, uv,v^2)$ and let $X$ denote the image of $\pi$.   Then, $\pi \colon \A^2 \to X$ is the normalization of $X$.  Suppose $f \in B$ agrees on the fibers of $\pi$. Write $f = \sum_{i=0}^m g_i(u)v^i, g_i \in K[u]$. Then $f(0,c) = f(0,-c)$ for all $c \in K$ and hence
$\sum_{i =0}^m g_i(0) v^i  = \sum_{i =0}^m (-1)^i g_i(0) v^i$.  Thus $g_i(0) = 0$ whenever $i$ is odd. Then
$f = \sum [g_i(u) - g_i(0)]v^i + \sum_{i \; \mathrm{even}} g_i(0)v^i$ is in $A$ since $g_i(u) - g_i(0) \in uB \subset A$.   This example generalizes to higher dimensions (see \cite[ Prop. 3.5]{Leahy:1981ta}.)  The zero set of $y^2 = x^2z$ is $X$ together with the negative $z$-axis and is called the Whitney umbrella.
Here is a sketch of the umbrella minus its handle,  which is the negative $z$-axis.
\index{Whitney umbrella}
\begin{figure}[h!]
\centering
\includegraphics[width=3in]{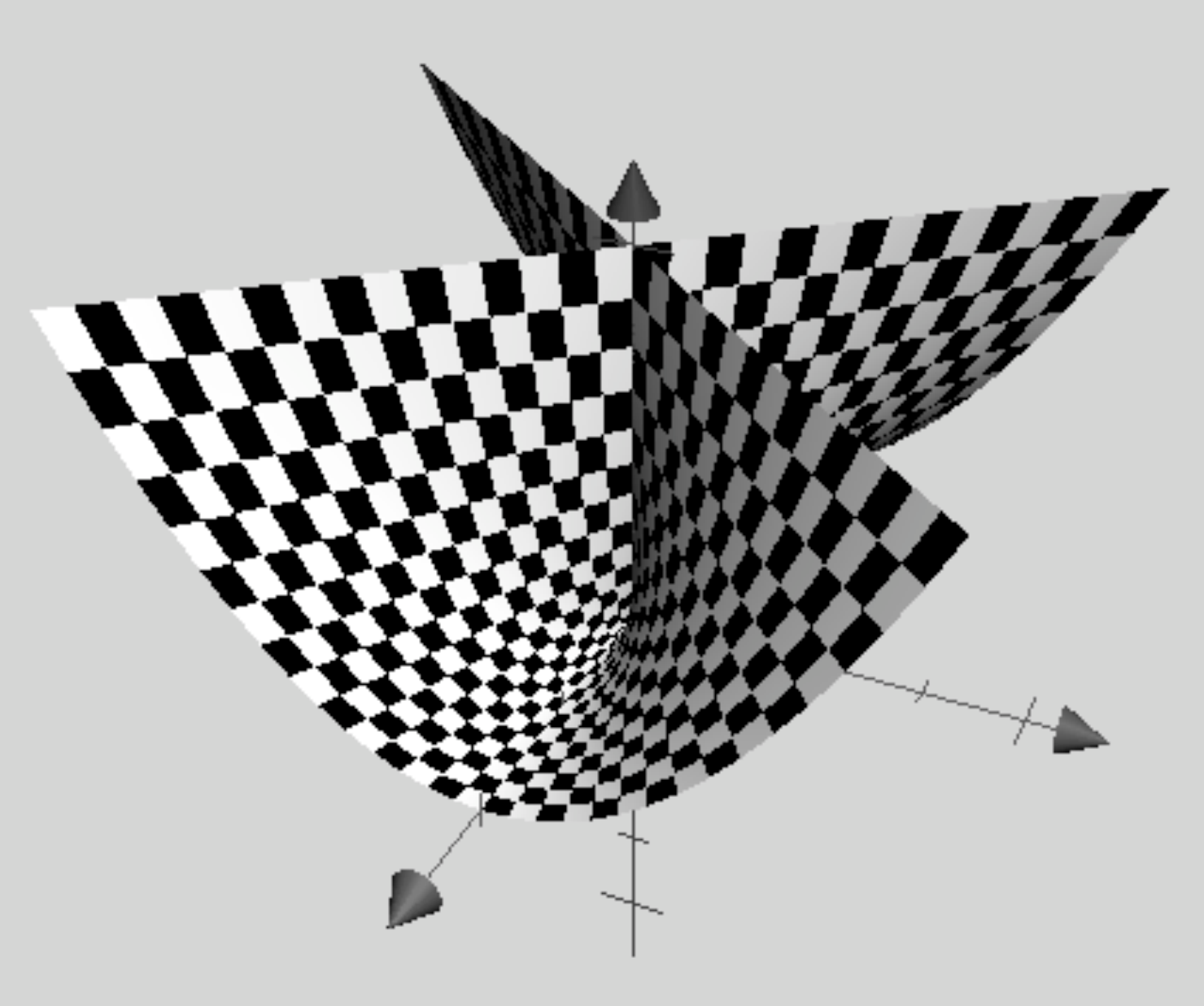}
 \caption{Whitney Umbrella}
\end{figure}

\end{ex}

Let $X$ be an algebraic variety over an algebraically closed field of characteristic 0. In the introduction we said  that a point $x \in X$  is a multicross if $x \in X$ is analytically isomorphic to $z \in Z$ where $Z$ is the union of linearly disjoint linear subspaces and $z$ is the origin.  We now recall the precise definition and cite a theorem, which asserts that a variety is weakly normal at a multicross singularity.   Then we cite a result about generic hyperplane sections of weakly normal varieties whose proof relied on various characterizations of the multicrosses due to Leahy and the current author.
\index{Leahy, John V.!multicross}  \index{Vitulli, Marie A.!multicross}

\begin{defn}[see \cite[Definition  3.3]{Leahy:1981ta} and \cite[Defintions 3.1]{Leahy:1981kr}] Let  $\mathcal{C}  = \{T_1, \ldots, T_r \}$ be a nonempty collection of disjoint subsets of $\{ 1, \ldots , p \}$ and let $z_1, \ldots, z_p$ be indeterminates.  We let $R_{\mathcal{C}}$ denote the complete local ring defined by
\begin{equation} 
R_{\mathcal{C}} = k[[z_1, \ldots , z_p]]/(z_{\alpha}z_{\beta} \mid \alpha \in T_i, \beta \in T_j, i \ne j \}.
\end{equation}
We say a  point $x$ on a variety $X$ is a multicross if $\widehat{\mathcal{O}}_{X,x}$ is isomorphic (as a $K$-algebra) to $R_{\mathcal{C}}$ for some $\mathcal{C}$ as above.
\end{defn}
\index{multicross}
We cite two fundamental results about the multicrosses.

\begin{prop}\cite[Prop. 3.4]{Leahy:1981ta} Let $X$ be an algebraic variety over an algebraically closed field of characteristic 0.  If $x \in X$ is a multicross then $\mathcal{O}_{X,x}$ is weakly normal.
\end{prop}

\begin{thm}\cite[Theorem 3.8]{Leahy:1981kr}  Suppose $X$ is a weakly normal variety over an algebraically closed field of characteristic 0 and let $Z$ denote the complement of the set of multicrosses.  Then, $Z$ is a closed subvariety of codimension at least two.
\end{thm}

The following result of the current author parallels  well-known results of Seidenberg \cite[see Theorems 7, $7'$, and 14]{MR0037548} on generic hyperplane sections of irreducible normal varieties.  The theorem was independently proven by C. Cumino, S. Greco, and M. Manaresi in \cite{MR700745}.  The conclusion of this theorem fails in positive characteristic by a class of examples that Cumino, Greco, and  Manaresi introduced in \cite{MR953739}.  In the statement of the next theorem, for a point $a \in \mathbb{A}^{m+1}$ by $H_a$ we mean the hyperplane $a_0 + a_1X_1 + \cdots + a_mX_m = 0$ in $\mathbb{A}^{m+1}$ .

\begin{thm} \cite[Theorem 3.4]{Vitulli:1983cr} Let $X \subset \mathbb{A}^m$ be an equidimensional weakly normal affine variety over an algebraically closed field of characteristic 0.  Then there exists a dense open subset $U$ of $\mathbb{A}^{m+1}$ such that $X \cap H_a$ is weakly normal whenever $a \in U$.
\end{thm}
\index{Vitulli, Marie A.!weak normality of hyperplane sections}

Instead of asking about whether a property that a variety has is also enjoyed by a general hyperplane section, one can ask that if a general hyperplane section has a property, does that property ``lift" to the original variety. 
Cumino and Manaresi defined a WN1 variety as a weakly normal variety such that the normalization $\overline{X} \to X$ is unramified in codimension one.  In \cite{MR953739} Cumino, Greco, and Manaresi show that if the general hyperplane sections are WN1 then so is the original variety. They used this fact to show that in positive characteristic, the general hyperplane section of a weakly normal variety is not weakly normal. 
\index{Cumino, Caterina!weak normality of hyperplane sections}  \index{Manaresi, Mirella!weak normality of hyperplane sections} \index{Greco, Silvio!weak normality of hyperplane sections} 
 \index{unramified}
 
 More recently R. Heitmann  proved a strong lifting result for \emph{seminormality}, which we now state.  Since we are dealing with general rings now we can not blur the distinction between seminormal and weakly normal rings.
 \index{Heitmann, Raymond!lifting property of seminormality}
 
\begin{thm}\cite[Main Theorem]{Heitmann:2008aa}
If $(R,\m)$ is a Noetherian local ring, $y$ is a regular element in $\m$, and $R/yR$ is seminormal, then $R$ is seminormal.
\end{thm}


\subsection{Weak Normality and the Chinese Remainder Theorem}\label{subsect:wn and CRT}

In Section \ref{subsec:CRT} we reported on seminormality and the Chinese Remainder Theorem.  We now mention a couple of results by Leahy and Vitulli that connect weak normality and the Chinese Remainder Theorem.  All varieties are taken over an algebraic closed field of characteristic 0.  The reader should see what these results say about affine rings and contrast the results to Theorems \ref{thm on CRT and sn}  and \ref{big thm on CRT and sn}. 
 \index{Leahy, John V.!CRT}  \index{Vitulli, Marie A.!CRT}
 
 For a closed subvariety $Y$ of an algebraic variety $X$ we let $\mathcal{I}_Y$ denote the sheaf of ideals defining $Y$.  For a point $X \in X$ we let $T_{X,x}$ denote the tangent space of $X$ at $x$.
 
 \begin{thm} \cite[Prop. 2.19]{Leahy:1981ta}  Let $X =X_1 \cup \cdots \cup X_n$ where each $X_i$ is a closed subvariety and suppose that $X_i$ is weakly normal for each $i$.  Further assume that $X_i \cap X_j = Y$ whenever $i \ne j$. Then $X$ is weakly normal if and only if $\mathcal{I}_Y = \mathcal{I}_{X_i} +  \mathcal{I}_{X_1 \cup \cdots \cup X_{i-1}}$ for $i = 1, \ldots, n$. 
 
 \end{thm}

If, in addition, we assume that the common intersection $X_i \cap X_j = Y$ is weakly normal we get the following result. This result can be used to show that a variety is weakly normal at a multicross singularity.

\begin{thm} Let $X = X_1 \cup \cdots \cup X_n$ be a union of closed subvarieties and assume that $X_i \cap X_j = Y$ whenever $i \ne j$, where $Y$ is weakly normal.  The following assertions are true.
\begin{enumerate}
\item $X$ is weakly normal if and only if each $X_i$ is weakly normal and $\mathcal{I}_Y = \mathcal{I}_{X_i} +  \mathcal{I}_{X_1 \cup \cdots \cup X_{i-1}}$ for $i = 1, \ldots, n$.
\item Suppose in addition that $Y$ is nonsingular.  Then $X$ is weakly normal if and only if each $X_i$ is weakly normal and $T_{Y,x} = T_{X_i,x} \cap T_{X_1 \cup \cdots \cup X_{i-1},x}$ for all $x \in Y$ and $i = 2, \ldots n$.

\end{enumerate}

\end{thm}

\subsection{The Weak Subintegral Closure of an Ideal} \label{subsec: wsi of ideal}

In this section we discuss the weak subintegral closure of an ideal. We use the definition proposed by the current author and Leahy \cite{Vitulli:1999fi}, which in turn is based on the 
criterion of Reid, Roberts and Singh \cite{Reid:1996dw}.  Our definition stands in the same relation to the definition of Reid, Roberts and Singh, as
the definition of the integral closure of an ideal does to the normalization of a ring.  In our definition we absorb the factor of $(-1)^n$ that appears in the equations (\ref{wsi eqns}).

\begin{defn}
Consider an $A$-ideal and a ring extension $I \subset A \subset B$. We say $b \in B$ is
{\bf weakly subintegral over
$I$}  provided that there exist
 $q \in \N$ and elements  $a_i\in{I}^i \;  (1\le{i}\le 2q+1$),
such that

\begin{equation}\label{LVeqns}
b^n+\sum_{i=1}^{n} {n \choose i} a_ib^{n-i}=0\quad (q+1\le n \le 2q+1).
\end{equation}

\noindent We let
$${}^*_BI = \{ b \in B \; | \; b \hskip 4pt{\rm is\hskip 4pt weakly\hskip 4pt
subintegral\hskip 4pt over\hskip 4pt }I\}.$$
We call ${}^*_BI  $ the {\bf weak subintegral closure of
$I$ in $B$}.  We write ${}^*I$ instead of ${}^*_AI  $ and refer to ${}^*I$ as the {\bf weak subintegral closure of $I$}. 
\end{defn}
\index{weakly subintegral over an ideal}
\index{${}^*_BI  $: weak subintegral closure of $I$ in $B$} 
\index{${}^*I  $: weak subintegral closure of $I$}
\index{Leahy, John V.!weak subintegral closure of ideal} \index{Vitulli, Marie A.!weak subintegral closure of ideal}

 The paper \cite{Vitulli:1999fi} contains an important link between weak normalization of a graded ring and weak subintegral closure of an ideal, which we recall. Suppose that $A \subseteq B$ are rings, $I$ is an ideal in $A$, and $b \in B$. Then, $b$ is weakly subintegral over $I^m$ if and only if the element $bt^m \in B[t]$ is weakly subintegral over the Rees ring $A[It]$ by  \cite[Lemma 3.2]{Vitulli:1999fi}.  Thus ${}^*_BI  $ is an ideal of ${}^*_BA$  \cite[Prop. 2.11]{Vitulli:1999fi}.  In particular, ${}^*I$ is an ideal of $A$.  Vitulli and Leahy also  showed that for an ideal $I$ in a reduced ring $A$ with finitely many minimal primes and total quotient ring $Q$, we have ${}^*(A[It]) = \oplus_{n \ge 0} \, {}^*_Q (I^n)t^n$ by \cite[Corollary 3.5]{Vitulli:1999fi}.
\index{Rees algebra}

Let's make a quick observation about a sufficient condition for an element to be weakly  normal over an ideal. 

\begin{lemma} \label{lemma high powers} Suppose that $I$ is an ideal of a ring $A$, $b \in A$ and $b^n \in I^n$ for all sufficiently high powers of $n$.  Then, $b \in {}^*I$.
\end{lemma}

\noindent \textbf{Proof.}
Suppose $q$ is such that $b^n \in I^n$ for all $n > q$.  Set $a_i = 0$ for $i= 1, \ldots, q$. 
Define $a_{q+1} = -b^{q+1}$.  Suppose $a_{q+1}, \ldots, a_{n-1}$ have been defined for $q+1 \le n -1 < 2q+1$.  
Set $a_n = -[b^n + \sum_{i=q+1}^{n-1} {n \choose i} a_ib^{n-i}]$. 
Since $a_n$ is an integer multiple of $b^n$, $a_n \in I^n$ and  $b^n+\sum_{i=1}^{n} {n \choose i} a_ib^{n-i}=0$.  By induction, we can define coefficients $a_i$ so that the equations (\ref{LVeqns}) are satisfied.  \qedbox

We can compare the weak subintegral closure of an ideal to what is usually called the Ratliff-Rush closure and was introduced by Ratliff and Rush in \cite{Ratliff:1978oa}.
\index{Ratliff-Rush closure}   \index{Rush, David E.}  \index{Ratliff, Louis J. Jr.}

\begin{cor}  Let $I$ be an ideal in a Noetherian ring $A$ containing a regular element of $A$ and let $\tilde{I}$ denote the Ratliff-Rush ideal associated with $I$, that is,  $\tilde{I} = \cup_{n \ge 0} (I^{n+1}:I^n)$.  Then, $\tilde{I} \subset {}^*I$.
\end{cor}

\noindent \textbf{Proof.}
This following immediately from the preceding corollary and the fact that $\tilde{I}^n = I^n$ for all sufficiently high powers of $n$, which was proven in \cite{Ratliff:1978oa}.  
\qedbox

\medskip

For an ideal $I$ of a ring $A$ we have inclusions $I \subset {}^*I \subset \overline{I} \subset \sqrt{I}$ and if $I$ is a regular ideal (i.e, $I$ contains a regular element) of a Noetherian ring we also have $I \subset \tilde{I} \subset  {}^*I \subset \overline{I} \subset \sqrt{I}$.
\index{regular ideal}

Vitulli and Reid \cite{Reid:1999sj} algebraically characterized weakly normal monomial ideals in a polynomial ring over a field.   Recall that the integral closure of a monomial ideal  $I$ in a polynomial ring $K[x_1,
\ldots, x_n]$ is generated by all monomials $x^{\gamma}$ such that $x^{m\gamma} \in I^m$  \emph{for some} positive integer $m$.   This is independent of the characteristic of $K$. The condition to be in the weak subintegral closure of a monomial ideal is slightly stronger and is characteristic dependent as we now explain.
\index{Reid, Les!weak subintegral closure of ideal}  \index{Vitulli, Marie A.!weak subintegral closure of monomial ideal}

\begin{prop} \cite[Proposition 3.3]{Reid:1999sj}  Let $I$ be a monomial ideal in a polynomial ring $K[x_1,
\ldots, x_n]$ in $n$ indeterminates over a field $K$. 
\begin{enumerate}
\item If char$(K) = 0$, then ${}^*I$ is the monomial ideal generated
by all monomials
$x^{\gamma}$ such that $x^{m\gamma} \in I^m$ for all sufficiently large positive integers
$m$.
\item  If char$(K)=p > 0$, then ${}^*I$ is the monomial ideal
generated by all monomials $x^{\gamma}$ such that $x^{p^m\gamma} \in I^{p^m}$
for some nonnegative integer $m$.
\end{enumerate}
\end{prop}

Reid and Vitulli \cite[Theorem 4.10]{Reid:1999sj} also presented a geometric characterization of the weak subintegral closure of a monomial ideal over a field of characteristic 0 in terms of the Newton polyhedron 
$\mathrm{conv}(\Gamma)$ of the exponent set $\Gamma = \Gamma(I)$ of the monomial ideal $I$.  We
 recall that $\Gamma(I)$ consists of all exponents of monomials in $I$ and $\mathrm{conv}(\Gamma)$ is
  the convex hull of $\Gamma = \Gamma(I)$.  We write $\overline{\Gamma}$ for the set of integral points
   in $\mathrm{conv}(\Gamma)$; thus $\overline{\Gamma}$ is the exponent set of the integral closure 
   $\overline{I}$ of $I$. We present a streamlined version of their characterization due to the current author.  First we recall the pertinent definitions.

\index{Newton polyhedron}  \index{$\mathrm{conv}(\Gamma)$: Newton polyhedron of $\Gamma(I)$}
\index{$\Gamma(I)$: exponent set of monomial ideal $I$}  \index{convex hull}
\index{$\overline{\Gamma}$}  
\index{$\overline{I}$: integral closure of ideal $I$}

\begin{defn} For a  polyhedron $P$ we define the {\bf relative interior} of $P$ by
$$\relint(P) := P - \cup E,$$ 
where the union is taken over all facets $E$ of $P$.
\end{defn}  
\index{$\relint(P)$: relative interior of polyhedron $P$} \index{relative interior of polyhedron}
\index{facet}

\begin{defn}  For a face $F$ of the Newton polyhedron $\Sigma=\conv(\Gamma)$  of a
monomial ideal, define
$${}^*F = \left\{ \mathbf{x} \in \relint(F) \hspace{.2cm} | \hspace{.2cm} \mathbf{x} = \sum n_i\gamma_i,
\hspace{.1cm}\gamma_i \in F \cap \Gamma \hspace{.1cm}, n_i
\in \mathbb{Z} \right\}.$$
That is, ${}^*F$ is the intersection of the group  generated by $F \cap
\Gamma$  with the relative interior of $F$.
\end{defn}

\begin{remark}  In \cite{Reid:1999sj}  the additional requirement that $\sum 
n_{i}=1$ was part of the definition of ${}^{*}F$.  Lemma 3.5 of \cite{Vitulli:2000jb} 
shows that this condition may and shall be deleted.   
\end{remark}
\index{${}^{*}F$}

\medskip
Before going further, we look at some examples.

\begin{ex} \label{ex 1 of *F}  Consider $I = (x^n, x^2y^{n-2}, y^n) \subset K[x,y]$, where $n = 2m+1$, let $\Gamma = \Gamma(I)$, and $\Sigma = \conv(\Gamma)$.   Notice that
\begin{eqnarray}
\label{ex eqn1} (n,0) &= &(0,n) +n(1,-1) \\
\label{ex eqn2} (2,n-2) &=& (0,n) +    2(1,-1).
\end{eqnarray}
Subtracting $m$ times equation (\ref{ex eqn2}) from equation (\ref{ex eqn1}) we get

\begin{equation}
(n,0)+(m-1)(0,n)-m(2,n-2)=(1,-1).
\end{equation}

With this one can check that every lattice point on the line segment $E$ from $(0,n)$ to $(n,0)$, which is a face of $\Sigma$, is in the group generated by $E \cap \Gamma$.   We may conclude that  ${}^*E = \{ (1, n-1), (2, n-2) , \ldots , ((n - 1, 1) \}$.
\end{ex}

\begin{ex} \label{ex 2 of *F}
Consider $I = (x^n, x^2y^{n-2}, y^n) \subset K[x,y]$, where $n$ is even,  let $\Gamma = \Gamma(I)$, and $\Sigma = \conv(\Gamma)$.   

Using the preceding example (after dividing all coordinates by 2)  we see that every lattice point with \emph{even} coordinates on the line segment $E$ from $(0,n)$ to $(n,0)$, which is a face of $\Sigma$, is in the group generated by $E \cap \Gamma$ so that ${}^*E = \{(2, n-2), \ldots , (n-2,2)\}$.
\end{ex}

Here is the geometric characterization of the weak subintegral closure of a monomial ideal in characteristic 0.

\begin{theorem} \label{geo char} Let $\Gamma = \Gamma(I)$ be the exponent set of a monomial ideal
$I$ in a polynomial ring $K[x_1, \dots, x_n]$ over a field $K$ of characteristic
0.  Then, an integral point $\gamma \in \conv(\Gamma)$ is in the exponent set
of ${}^*I$ if and only if $\gamma \in \bigcup {}^*F$, where the union is taken
over all faces $F$ of $\conv(\Gamma)$.
\end{theorem}

We include an example that illustrates the distinction between the integral closure and weak subintegral closure of a bivariate monomial ideal over a field $K$ of characteristic 0.  By pos$(X)$  we mean the positive cone of the subset $X$ of $\mathbb{R}^n$.
\index{$\mathrm{pos}(X)$: positive cone of $X$}  
 \index{positive cone of $X$}

\index{${}^*\Gamma$}

\begin{ex} \label{wn of mono ideal ex}  Let $I = (x^6, x^2y^4, y^6) \subset K[x,y]$, $\Gamma=\Gamma(I)$,  ${}^*\Gamma=\Gamma({}^*I)$ and $\Sigma
= \mathrm{conv}(\Gamma)$.  The exponent set $\Gamma$ consists of all lattice points
on or above the thick-lined staircase figure and $\Sigma$ is the 1st quadrant with the lower left-hand
corner clipped, i.e.,  $\Sigma$ consists of all points on or above the oblique line joining $V_1$ and $V_2$.  Observe that $\Sigma$ and is the sole 2-dimensional face of
$\mathrm{conv}(\Gamma)$,  $E_1= \{ 0\} \times [6, \infty )  $, $E_2=[(0,6),(6,0)]$, and
$E_3= [6, \infty ) \times \{0 \} $ are the edges  and
$V_1=\{(0,6)\}$ and $V_2=\{(6,0)\}$ are the vertices of $\Sigma$.  The sets $\Gamma$,
$\Sigma$, $E_1, E_2, E_3, V_1, V_2$ are depicted below. The lattice points depicted by open circles are in $\Sigma \setminus {}^*\Gamma$.  The lattice points depicted by filled circles are in ${}^*\Gamma$.
 
 \begin{center}
\setlength{\unitlength}{5mm}
\begin{picture}(10,10)
\linethickness{.1mm}
\multiput(0,0)(0,1){8}{\line(1,0){8}}
\multiput(0,0)(1,0){8}{\line(0,1){8}}
\linethickness{.5mm}
\put(0,6){\line(0,1){2}}
\put(6,0){\line(1,0){2}}
\put(0,6){\line(1,0){2}}
\put(2,6){\line(0,-1){2}}
\put(2,4){\line(1,0){4}}
\put(6,0){\line(0,1){4}}
\put(6,0){\line(-1,1){6}}
\put(-1,6){$V_1$}
\put(0,6){\circle*{.3}}
\put(6,0){\circle*{.3}}
\put(2,4){\circle*{.3}}
\put(4,2){\circle*{.3}}
\put(4,3){\circle*{.3}}
\put(5,2){\circle*{.3}}
\put(5,3){\circle*{.3}}
\put(1,5){\circle{.3}}
\put(3,3){\circle{.3}}
\put(5,1){\circle{.3}}
\put(.5,1){\mbox{ \huge $\Sigma$}}
\put(0,-2.5){Fig. 3:  Weak Subintegral Closure of $(x^6, x^2y^4,y^6)$}
\put(4.0,5.0){\mbox{ \huge $\Gamma$}}
\put(-1.5,7){\mbox{ \large $E_1$}}
\put(1.5,2,5){\mbox{ \large $E_2$}}             
\put(7,-1.5){\mbox{ \large $E_3$}}
\put(6,-1){$V_2$}
\end{picture}
\end{center}
\bigskip
\bigskip
\medskip
\bigskip
Observe that
\begin{eqnarray*}
 {}^*\Sigma &=& \left( \Gamma \setminus (E_1 \cup E_2 \cup E_3) \cap \Gamma \right) \cup
\{ (4,3),(5,2), (5,3) \}; \\
 {}^*E_1 &=& \{ (0,7),(0,8),(0,9), \dots \}; \\
 {}^*E_2 &=& \{ (2,4),(4,2) \}; \\
 {}^*E_3 &=& \{ (7,0),(8,0),(9,0), \dots \}; \\
 {}^*V_i &=& V_i \quad (i=1,2); \\
 {}^*\Gamma &=& \Gamma \cup \{ (4,2),(4,3),(5,2),(5,3) \}; \mbox{  and } \\
\overline{\Gamma} &=& \Gamma \cup \{ (1,5),(3,3),(4,2),(4,3),(5,1),(5,2),(5,3) \}.
\end{eqnarray*}
Thus we have $\Gamma \subset {}^*\Gamma \subset \overline{\Gamma}$, where both containments are proper.
\end{ex}

Notice that all of the points in the relative interior of the Newton polyhedron are in the weak subintegral closure of $I$.  An analog of this observation is true for arbitrary ideals (see Prop. \ref{prop:lantz} below), which we will present after we introduce some necessary notation.  For more details the reader should consult \cite{GaffVit2008:aa}.

Let $I$ be a  monomial ideal over a field of arbitrary characteristic.  The proof of  Theorem \ref{geo char} shows that if  $\alpha$ is in ${}^*F$ for some face $F$ of $\conv(\Gamma(I))$, then   $x^{n\alpha} \in I^n$ for all sufficiently large $n$.  Thus we may conclude that $x^{\alpha} \in {}^*I$ by Lemma \ref{lemma high powers}. We now  take a look at the ideal with the same monomial generators as in Example \ref{wn of mono ideal ex} but over a field of characteristic 2.  

\begin{ex}  Assume that $\mathrm{char}(K)=2$, let $I = (x^6, x^2y^4, y^6) \subset K[x,y]$, $\Gamma=\Gamma(I)$,  ${}^*\Gamma=\Gamma({}^*I)$ and $\Sigma   = \mathrm{conv}(\Gamma)$. Notice that
\begin{eqnarray*}
\left( xy^5 \right)^2 &=& \left( x^2y^4 \right) y^6 \in I^2  \\
\left( x^3y^3 \right)^2 &=& x^6y^6 \in I^2   \\
\left( xy^5 \right)^2 &=&\left(x^2y^4 \right)y^6 \in I^2
\end{eqnarray*}
Combining these calculations with those in Example \ref{ex 2 of *F} we see that ${}^*I = \overline{I}$.  This explicit example illustrates that unlike integral closure,  the weak subintegral closure of a  monomial ideal depends on the characteristic of the base field.
\end{ex}

We now turn our attention to more general ideals in Noetherian rings.

\begin{notn}  For an ideal $I$ of a Noetherian ring $A$ and an
element $a \in A$ we write
$\ord_I(a) = n$ if $a \in I^n \setminus I^{n+1}$ and $\ord_I(a)
= \infty$ if $a \in \bigcap_{n \ge 1} I^n$.  Next we define 
$$\overline{v}_I(a) = \lim_{n \to \infty} \frac{\ord_I(a^n)}{n}.$$  
The indicated limit always exists (possibly
being $\infty$; (\cite[Prop. 11.1]{MR722609} or \cite[Theorem 10.1.6]{MR2266432}) and 
$\overline{v}_I$ is called the \textbf{asymptotic Samuel function }of $I$.  
\end{notn}
\index{asymptotic Samuel function} \index{$\overline{v}_I$: asymptotic Samuel function of $I$}

\medskip

 Let $I$ be a regular  ideal in a Noetherian ring $A$  and $a \in A$.  The asymptotic Samuel function $\overline{v}_I$ is determined by the Rees valuations $v_j$ of $I$.  Namely,
$$\overline{v}_I(a) = \min_j \left\{ \frac{v_j(a)}{v_j(I)}
\right\},$$ where $v_j(I) = \min \{ v_j(b) \mid b \in I \} $ (see \cite[Lemma 10.1.5]{MR2266432}). We refer the reader to Chapter 10 of  \cite{MR2266432} for the fundamentals on Rees valuations of ideals. Recall that $\overline{v}_I = \overline{v}_J$ whenever $\overline{J}  = \overline{I}$ (see \cite[Cor. 11.9]{MR722609}).  This immediately implies that $J_> = I_>$ whenever $\overline{J}  = \overline{I}$.
\medskip
\index{Rees valuations}

\begin{notn}
For an ideal $I$ in a Noetherian ring $A$ we let
$$I_> =  \{ a \in A \mid \overline{v}_I(a) > 1 \}.$$  
\end{notn}
\index{$I_>$: ideal associated with $I$}

By elementary properties of the asymptotic Samuel function,  $I_>$ is an ideal of $A$ and 
a subideal of $\overline{I}$.   It need not contain the original ideal $I$ as we now illustrate.  

\begin{ex}  Let $K$ be a field and $I = (x^6, x^2y^4, y^6) \subset K[x,y]$.  Then, $I$ has one Rees valuation, namely the monomial valuation defined by $v(x^ay^b)=a+b$ (see Chapter 10 of  \cite{MR2266432} for the fundamentals on monomial valuations).  We have $v(I) = 6$ and $I_> = \{f \in K[x,y] \mid v(f) > 6 \} = (x,y)^7$. In particular, $I_>$ does not contain $I$.

\end{ex}

The ideal $I_>$ plays an important role in conditions from stratification theory such as Whitney's condition A and Thom's condition $A_f$; the reader can learn more about these conditions in \cite{MR1702106}. To give the reader a little more feeling for the ideal $I_>$ we cite some lemmas that were proven in  \cite{GaffVit2008:aa}.

\begin{lemma}\cite[Lemma 4.2 ]{GaffVit2008:aa}   Let $I$ be a regular ideal in a Noetherian ring $A$. 
Then,
$$I_> = \bigcap_i  \m_iIV_i \cap A,$$
where the intersection is taken over all Rees valuation rings $(V_i, \m_i)$ of $I$.
In particular, $I_>$ is an integrally closed ideal.
\end{lemma}
\index{regular ideal}

\begin{lemma}\cite[Lemma 4.3 ]{GaffVit2008:aa}  Let $I$ be a nonzero monomial ideal in a polynomial ring over a field. Then, $I_>$ is again a monomial ideal.
\end{lemma}

We now present a generalization of what was known as Lantz's conjecture (after a talk D. Lantz gave in 1999 at a Route 81 Conference in upstate NY) and illustrate the result with an example. Lantz conjectured that if $I$ is an $\m$-primary ideal in  a 2-dimensional regular local ring $(A, \m)$, then ${}^*I$ contains all elements $a \in A$ such that $\overline{v}_I(a)  > 1$.
\index{Lantz, David!Lantz's conjecture}
\index{Vitulli, Marie A.!Lantz's conjecture}  \index{Gaffeny, Terence!Lantz's conjecture}

 \begin{prop}\label{prop:lantz} \cite[Prop. 4.4]{GaffVit2008:aa}  Let $I$ be an ideal of a
Noetherian  ring $A$.  Then, $I_> \subseteq {}^*I$.
\end{prop}

\begin{ex}  Let $K$ be a field of  characteristic 0, let $n=2m+1$, and consider $I = (x^n, x^2y^{n-2}, y^n)$.  We claim that ${}^* I = \overline{I}$.   By Example \ref{ex 1 of *F}  we may  conclude that for the facet $E$ of conv$(\Gamma)$, ${}^*E = \{(1,n-1), (2,n-2), \ldots, (n-2,2), (n-1,1) \}$. This, together with the fact that $I_> \subset {}^*I$, implies that ${}^* I = \overline{I}$.
\end{ex}

\begin{ex}  Let $K$ be a field of characteristic 0, let $n=2m$, and consider $I = (x^n, x^2y^{n-2}, y^n)$.
Using Example \ref{ex 2 of *F}  and Theorem \ref{geo char}  we see that the weak subintegral closure is ${}^*I = (x^n, x^{n-2}y^2, \ldots , x^2y^{n-2}, y^n )$.
\end{ex}

One consequence of the preceding proposition is the following connection between the ideal $I_>$ and the minimal reductions $\mathcal{MR}(I)$ of the ideal $I$.
\index{$\mathcal{MR}(I)$: minimal reductions of $I$}

\begin{cor}  \label{cor:lantz}  \cite[Cor. 4.5 ]{GaffVit2008:aa} Let $I$ be an ideal of a Noetherian  ring $A$. 
Then, $$I_> \subseteq \bigcap_{J \in \mathcal{MR}(I)} {}^*J.$$
\end{cor}

We mention another interesting occurrence of the ideal $I_>$.  For more results in the same spirit the reader should consult \cite{GaffVit2008:aa}.

\begin{theorem}  \label{thm:weaknor of min red} \cite[Theorem 4.6 ]{GaffVit2008:aa} Let $(A, \mathfrak{m},k)$ be a Noetherian local
ring  such that $k$ is algebraically closed of
characteristic 0.  Suppose that $I$ is an
$\mathfrak{m}$-primary ideal.   If $J$ is any minimal reduction
of $I$ then $J + I_> = {}^*J$.
\end{theorem}

It is well known that the integral closure of an ideal $I$ in an integral domain $A$ can be characterized in terms of valuation rings.  Namely, $\overline{I} = \bigcap _V IV \cap A$, where the intersection is taken over all valuation rings of the quotient field of $A$ (see \cite[Prop. 6.8.1]{MR2266432}).    This can be rephrased in terms of maps into valuation rings and for Noetherian rings in terms of maps into discrete rank one valuation rings.  More precisely,   an element $f$ in a Noetherian ring $A$ is in $\overline{I}$ if and only if $\rho(f) \in \rho(I)V$ for every homomorphism $ \rho \colon A \to V$, where $V$ is a  discrete rank one valuation ring.  For algebro-geometric local rings that we can limit which discrete rank one valuation rings we look at. The following is an analog of the complex analytic criterion involving germs of morphisms from the  germ of the pointed complex unit disk $(\mathbb{D}, 0)$ to the  germ $(X,x)$ (e.g., see \cite[2.1 Theor\`eme]{LejeuneTeiss:1974bh}).

\begin{prop} \label{alg curve crit} \cite[Prop. 5.4 ]{GaffVit2008:aa}  Let $(A, \m, k)$ be the local ring of an algebraic variety over an algebraically closed field $k$, $I$ an ideal of $A$, and $h \in A$.  Then, $h \in \overline{I} \Leftrightarrow$ for every local homomorphism of $k$-algebras $\rho \colon A \to k[[X]]$ we have 
$\rho(h) \in \rho(I)k[[X]]$. 
\end{prop} 

In the criteria involving maps into discrete valuation rings, one thinks of the targets of those maps as test rings for determining integral dependence.  Recently H. Brenner \cite{Brenner2006aa} suggested a valuative criterion for when an element is in the weak subintegral closure of an ideal in an algebro-geometric local ring using certain 1-dimensional seminormal local rings as test rings and the results of \cite{GaffVit2008:aa}.  Gaffney-Vitulli  took another approach to developing a valuative criterion in both the algebraic and complex analytic settings (see \cite[Prop. 5.7 and Prop. 5.8]{GaffVit2008:aa}).  
\index{Brenner, Holger!valuative criterion}  \index{Gaffney, Terence!valuative criterion} \index{Vitulli, Marie A.!valuative criterion}





We end our account here.  We hope we have left you, the reader, with a better idea of the many ramifications of the closely related notions of weak and seminormality.



\printindex

\end{document}